\documentclass[12pt,reqno]{amsart} 
\usepackage{amssymb}\usepackage{epsf}

\newtheorem{theorem}{Theorem}
\theoremstyle{definition}
\DeclareMathOperator*{\adj}{adj}
\DeclareMathOperator*{\ind}{ind}

\newcommand{\eg}{\Gamma}
\newcommand{\be}{\begin{equation}}
\newcommand{\ee}{\end{equation}}
\newcommand{\ba}{\begin{eqnarray}}
\newcommand{\ea}{\end{eqnarray}}

\newcommand{\lab}[1]{\label{#1}}

\newcommand{\C}{\mathbb C}
\newcommand{\R}{\mathbb R}
\newcommand{\Z}{\mathbb Z}

\newcommand{\T}{\mathbb T}

\newcommand{\ve}{\varepsilon}

\begin{document}

\title
{Elliptic hypergeometric functions}

\author{V. P. Spiridonov}

\address{Laboratory of theoretical physics, JINR,
Dubna, Moscow reg., 141980, Russia}

\thanks{
This is a complementary chapter to the book by
G. E. Andrews, R. Askey, and R. Roy,
{\em Special Functions}, Encyclopedia of Math. Appl.
{\bf 71}, Cambridge Univ. Press, Cambridge, 1999,
written for its Russian edition:
Moscow, MCCME, 2013, pp. 577-606.
}

\maketitle

{\bf \em Introduction.}
The wonderful book by Andrews, Askey, and Roy \cite{aar} is mainly devoted
to special functions of hy\-per\-geo\-met\-ric type -- to the plain hy\-per\-geo\-met\-ric
series and integrals and their $q$-analogues. Shortly before its publication
there appeared first examples of hypergeometric functions of a new type
related to elliptic curves.
A systematic theory of elliptic hy\-per\-geo\-met\-ric functions was constructed
in 2000-2004 over a short period of time. The present complement
reviews briefly the status of this theory by the spring of 2013. It repeats
where possible the structure of the book \cite{aar}, and it
is substantially based on author's thesis \cite{spi:thesis} and survey \cite{spi:umnrev}.

The theory of quantum and classical completely integrable systems played
a crucial role in the discovery of these new special functions.
An elliptic extension of the terminating very-well-poised balanced
$q$-hy\-per\-geo\-met\-ric series $_{10}\varphi_9$ with discrete values of
parameters appeared for the first time in elliptic solutions
of the Yang-Baxter equation \cite{ft} associated with the exactly
solvable models of statistical mechanics \cite{djkmo}.
The same terminating series with arbitrary parameters appeared in
\cite{spi-zhe:spectral} as a particular solution of a pair of
linear finite difference equations, the compatibility condition of which
yields the most general known $(1+1)$-dimensional nonlinear
integrable chain analogous to the discrete time Toda chain.
An elliptic analogue of Euler's gamma function
depending on two bases $p$ and $q$ of modulus less than 1,
which already appeared in Baxter's eight vertex model
\cite{bax}, was investigated in \cite{rui:first},
and in \cite{spi:theta2} a modified elliptic gamma function
was constructed for which one of the bases may lie on the unit circle.
General elliptic hypergeometric functions are defined by the integrals
discovered in \cite{spi:umn}, which qualitatively differ from the
terminating elliptic hypergeometric series.
The appearance of such mathematical objects was quite unexpected, since
no handbook or textbook of special functions contained any hint of their
existence. However, the generalized gamma functions
related to elliptic gamma functions and forming one of the key ingredients
of the theory were constructed long ago
by Barnes \cite{bar:multiple} and Jackson \cite{jac:basic}.
The most important known application of the elliptic hypergeometric
integrals was found quite recently -- they emerged in the description
of topological characteristics of four-dimensional supersymmetric
quantum field theories \cite{DO,gprr,SV,SV2}.

\medskip

{} {\bf \em Generalized gamma functions.} In the beginning of XXth
century Barnes \cite{bar:multiple} constructed the following multiple zeta function:
$$
\zeta_m(s,u;\mathbf{\omega})=\sum_{n_1,\ldots,n_m\in\Z_{\geq 0}} \frac{1}{(u+\Omega)^s},
\quad\Omega=n_1\omega_1+\cdots+n_m\omega_m, \quad \Z_{\geq 0}=0,1,\ldots,
$$
where $u,\ \omega_j\in\C$. This series converges for $\text{Re}(s)>m$ provided all
$\omega_j$ lie on one side of a line passing through the point $u=0$
(this forbids accumulation points of the $\Omega$-lattice in compact domains).
Using an integral representation for analytical continuation of $\zeta_m$ in $s$,
Barnes also defined the multiple gamma function
$\Gamma_m(u;\mathbf{\omega}) =\exp(\partial \zeta_m(s,u;\mathbf{\omega})/\partial
s)|_{s=0}.$ It has the infinite product representation
\begin{equation}
\frac{1}{\Gamma_m(u;\mathbf{\omega})} =
e^{\sum_{k=0}^m\gamma_{mk}\frac{u^k}{k!}}\;u\makebox[-1em]{}
\sideset{}{'}\prod_{n_1,\dots,
n_m\in\Z_{\geq 0}} \left(1+\frac{u}{\Omega}\right)
e^{\sum_{k=1}^m(-1)^k\frac{u^k}{k\Omega^k}}, \label{m-gamma}\end{equation} where
$\gamma_{mk}$ are some constants analogous to Euler's constant (in \cite{bar:multiple}, the
normalization $\gamma_{m0}=0$ was used). The primed product means that the point
$n_1=\ldots=n_m=0$ is excluded from it. The function $\Gamma_m(u;\mathbf{\omega})$
satisfies $m$ finite difference equations of the first order
\begin{equation}
\Gamma_{m-1}(u;\mathbf{\omega}(j))\Gamma_m(u+\omega_j;\mathbf{\omega})
=\Gamma_m(u;\mathbf{\omega}), \qquad j=1,\ldots,m, \label{bar-eq}\end{equation} where
$\mathbf{\omega}(j)=(\omega_1,\ldots,\omega_{j-1},\omega_{j+1},\ldots, \omega_m)$ and
$\Gamma_0(u;\omega):=1/u$. The function $\Gamma_1(\omega_1 x;\omega_1)$
essentially coincides with the Euler gamma function $\Gamma(x)$.
The plain, $q$-, and elliptic hy\-per\-geo\-met\-ric functions are
connected to $\Gamma_m(u;\mathbf{\omega})$ for $m=1,2,3$, respectively.

Take $m=3$ and assume that $\omega_{1,2,3}$ are pairwise
incommensurate quasiperiods. Then define three base variables:
\begin{eqnarray*}
&& q= e^{2\pi\textup{i}\frac{\omega_1}{\omega_2}}, \quad
p=e^{2\pi\textup{i}\frac{\omega_3}{\omega_2}}, \quad  r=e^{2\pi\textup{i}\frac{\omega_3}{\omega_1}},
\\ &&
\tilde q= e^{-2\pi\textup{i}\frac{\omega_2}{\omega_1}}, \quad
\tilde p=e^{-2\pi\textup{i}\frac{\omega_2}{\omega_3}},   \quad
\tilde r=e^{-2\pi\textup{i}\frac{\omega_1}{\omega_3}},
\end{eqnarray*}
where $\tilde q,\tilde p,\tilde r$ denote the $\tau\to -1/\tau$
modular transformed bases. For $|p|,|q|<1$, the infinite products
$$
(z;q)_\infty=\prod_{j=0}^\infty (1-zq^j), \qquad
(z;p,q)_\infty=\prod_{j,k=0}^\infty (1-zp^jq^k)
$$
are well defined. It is easy to derive equalities \cite{jac:basic}
\begin{equation}
\frac{(z;q)_\infty}{(qz;q)_\infty}=1-z, \qquad
\frac{(z;q,p)_\infty}{(qz;q,p)_\infty}=(z;p)_\infty, \qquad
\frac{(z;q,p)_\infty}{(pz;q,p)_\infty}=(z;q)_\infty.
\label{eq-1}\end{equation}

The odd Jacobi theta function (see formula (10.7.1)
in \cite{aar}) can be written as
\begin{align*}
\theta_1( u|\tau) &=
-\textup{i}\sum_{n=-\infty}^\infty (-1)^n e^{\pi\textup{i}\tau(n+1/2)^2}
e^{\pi\textup{i}(2n+1)u}
\\
&= \textup{i} p^{1/8} e^{-\pi\textup{i} u}\: (p;p)_\infty\: \theta(e^{2\pi\textup{i} u};p), \quad
u\in\mathbb{C},
\end{align*}
where $p=e^{2\pi\textup{i}\tau}$. The modified theta function (see Theorem 10.4.1
in \cite{aar})
\begin{equation}
\theta(z;p):=(z;p)_\infty (pz^{-1};p)_\infty
 =\frac{1}{(p;p)_\infty} \sum_{k\in\Z} (-1)^kp^{k(k-1)/2}z^k
\label{short-theta}\end{equation}
plays a crucial role in the following. It obeys the following properties:
\begin{equation}
\theta(pz;p)=\theta(z^{-1};p)=-z^{-1}\theta(z;p)
\label{theta-trafo}\end{equation}
and $\theta(z;p)=0$ for $z=p^k,\; k\in\Z.$ We denote
$$
\theta(a_1,\ldots,a_k;p):=\theta(a_1;p)\cdots\theta(a_k;p), \quad
\theta(at^{\pm 1};p):=\theta(at;p)\theta(at^{-1};p).
$$
Then the Riemann relation for products of four theta functions takes the form
\begin{eqnarray}
\theta(xw^{\pm 1},yz^{\pm 1};p) -\theta(xz^{\pm 1},yw^{\pm 1};p)
=yw^{-1}\theta(xy^{\pm 1},wz^{\pm 1};p)
\label{ident}\end{eqnarray}
(the ratio of the left- and  right-hand sides is a bounded function of
the variable $x\in\C^*$,
and it does not depend on $x$ due to the Liouville theorem,
but for $x=w$ the equality is evident).

Euler's gamma function can be defined as a special meromorphic
solution of the functional equation $f(u+\omega_1)=uf(u)$. Respectively, $q$-gamma functions
are connected to solutions of the equation
$f(u+\omega_1)=(1-e^{2\pi\textup{i} u/\omega_2})f(u)$ with $q=e^{2\pi\textup{i}\omega_1/\omega_2}$.
For $|q|<1$, one of the solutions has the form
$\Gamma_q(u)=1/(e^{2\pi\textup{i} u/\omega_2};q)_\infty$
defining the standard $q$-gamma function (it differs from function
(10.3.3) in \cite{aar} by the substitution $u=\omega_1 x$ and some
elementary multiplier). The modified $q$-gamma function (``the double sine",
``non-compact quantum dilogarithm",``hyperbolic gamma function"),
which remains well defined even for  $|q|=1$, has the form
\be
\gamma(u;\mathbf{\omega})=
\exp\left(-\int_{\R+\textup{i}0}\frac{e^{ux}}
{(1-e^{\omega_1 x})(1-e^{\omega_2 x})}\frac{dx}{x}\right),
\label{int-q-gamma}\ee
where the contour $\R+\textup{i}0$ coincides with the real axis deformed to pass clockwise
the point $x=0$ in an infinitesimal way. If
$\text{Re}(\omega_1), \text{Re}(\omega_2)>0$, then the integral converges for
$0<\text{Re}(u)< \text{Re}(\omega_1+\omega_2)$.
Under appropriate restrictions on $u$ and $\omega_{1,2}$, the
integral can be computed as a convergent sum of the residues of poles
in the upper half plane. When $\text{Im}(\omega_1/\omega_2)>0$, this yields
the expression
$\gamma(u;\mathbf{\omega})=(e^{2\pi\textup{i}u/\omega_1}\tilde q;\tilde q)_\infty
/ (e^{2\pi\textup{i} u/\omega_2}; q)_\infty$, which can be extended analytically
to the whole complex $u$-plane. This function, serving as a key building block of
the $q$-hy\-per\-geo\-met\-ric functions for $|q|=1$, was not considered in
 \cite{aar} and \cite{gas-rah:basic}; for its detailed description
see \cite{faddeev,jm,kls:unitary,rui:first,volk} and the
literature cited therein.

In an analogous manner, elliptic gamma functions are connected to
the equation
\be
f(u+\omega_1)=\theta(e^{2\pi\textup{i} u/\omega_2};p)f(u).
\lab{e-gamma-eq}\ee
Using the factorization \eqref{short-theta} and equalities \eqref{eq-1},
it is easy to see that the ratio
\begin{equation}
 \eg(z;p,q) = \frac{(pqz^{-1};p,q)_\infty}{(z;p,q)_\infty}
\label{ell-gamma}\end{equation}
satisfies the equations
$$
\eg(qz;p,q)=\theta(z;p)\eg(z;p,q),\quad \eg(pz;p,q)=\theta(z;q)\eg(z;p,q).
$$
Therefore the function $f(u)=\eg(e^{2\pi\textup{i} u/\omega_2};p,q)$ defines a
solution of equation \eqref{e-gamma-eq} valid for $|q|, |p|<1$, which is called
the (standard) elliptic gamma function \cite{rui:first}. It can be
defined uniquely as a meromorphic solution of three equations:
equation \eqref{e-gamma-eq} and
$$
f(u+\omega_2)=f(u),\qquad f(u+\omega_3)=\theta(e^{2\pi\textup{i} u/\omega_2};q)f(u)
$$
with the normalization $f(\sum_{m=1}^3\omega_m/2)=1$, since
non-trivial triply periodic functions do not exist. The reflection formula
has the form $\eg(z;p,q)\eg(pq/z;p,q)=1$.
For $p= 0$, we have $\eg(z;0,q)=1/(z;q)_\infty$.

The modified elliptic gamma function, which is well defined
for $|q|=1$, has the form \cite{spi:theta2}
\be
G(u;\mathbf{\omega})=
\eg(e^{2\pi\textup{i} \frac{u}{\omega_2}};p,q)
\Gamma(re^{-2\pi\textup{i} \frac{u}{\omega_1}};\tilde q,r).
\lab{unit-e-gamma}\ee
It yields the unique solution of three equations: equation \eqref{e-gamma-eq} and
$$
f(u+\omega_2) =\theta(e^{2\pi\textup{i} u/\omega_1};r) f(u),
\qquad   f(u+\omega_3) =e^{-\pi\textup{i}B_{2,2}(u;\mathbf{\omega})} f(u)
$$
with the normalization $f(\sum_{m=1}^3\omega_m/2)=1$. Here
$$
B_{2,2}(u;\mathbf{\omega})=\frac{u^2}{\omega_1\omega_2}
-\frac{u}{\omega_1}-\frac{u}{\omega_2}+
\frac{\omega_1}{6\omega_2}+\frac{\omega_2}{6\omega_1}+\frac{1}{2}
$$
denotes the second order Bernoulli polynomial appearing in the modular transformation
law for the theta function
\begin{equation}
\theta\left(e^{-2\pi\textup{i}\frac{u}{\omega_1}};
e^{-2\pi\textup{i}\frac{\omega_2}{\omega_1}}\right)
=e^{\pi\textup{i}B_{2,2}(u;\mathbf{\omega})} \theta\left(e^{2\pi\textup{i}\frac{u}{\omega_2}};
e^{2\pi\textup{i}\frac{\omega_1}{\omega_2}}\right).
\label{mod-theta}\end{equation}
One can check \cite{die-spi:unit} that the same three equations
and normalization are satisfied by the function
\be
G(u;\mathbf{\omega})
= e^{-\frac{\pi \textup{i}}{3}B_{3,3}(u;\mathbf{\omega})}
\Gamma(e^{-2\pi\textup{i} \frac{u}{\omega_3}};\tilde r,\tilde p),
\lab{mod-e-gamma}\ee
where $|\tilde p|,|\tilde r|<1$, and  $B_{3,3}(u;\mathbf{\omega})$
is the third order Bernoulli polynomial
$$
B_{3,3}\Big(u+\sum_{m=1}^3\frac{\omega_m}{2};\mathbf{\omega}\Big)
=\frac{u(u^2-\frac{1}{4}\sum_{m=1}^3\omega_m^2)}{\omega_1\omega_2\omega_3}.
$$
The functions \eqref{unit-e-gamma} and \eqref{mod-e-gamma} therefore coincide,
and their equality defines one of the laws of
the $SL(3,\Z)$-group of modular transformations for the elliptic gamma function
\cite{fel-var:elliptic}. From expression \eqref{mod-e-gamma}, the function
$G(u;\mathbf{\omega})$ is seen to remain meromorphic when
$\omega_1/\omega_2>0$, i.e. when $|q|=1$. The reflection formula
for it has the form
$G(a;{\bf \omega})G(b;{\bf \omega})=1,$ $a+b=\sum_{k=1}^3\omega_k.$
In the regime $|q|<1$ and $p,r\to0$ (i.e., Im$(\omega_3/\omega_1)$,
Im$(\omega_3/\omega_2)\to +\infty$), expression \eqref{unit-e-gamma}
obviously degenerates to the modified $q$-gamma function $\gamma(u;{\bf \omega})$.
Representation \eqref{mod-e-gamma} yields an alternative way of reduction
to $\gamma(u;{\bf \omega})$; a rigorous limiting connection of
such a type was built for the first time in a different way
by Ruijsenaars \cite{rui:first}.

As shown by Barnes, the $q$-gamma function
$1/(z;q)_\infty$ where $z=e^{2\pi\textup{i} u/\omega_2}$ and
$q=e^{2\pi\textup{i}\omega_1/\omega_2}$, $\text{Im} (\omega_1/\omega_2)>0$,
equals the product $\Gamma_2(u;\omega_1,\omega_2)\Gamma_2(u-\omega_2;\omega_1,-\omega_2)$
up to the exponential of a polynomial.
Similarly, the modified $q$-gamma function $\gamma(u;\mathbf{\omega})$
equals up to an exponential factor to the ratio
$\Gamma_2(\omega_1+\omega_2-u;\mathbf{\omega})/\Gamma_2(u;\mathbf{\omega})$.
Since $\theta(z;q)=(z;q)_\infty(qz^{-1};q)_\infty$, the
$\Gamma_2(u;\mathbf{\omega})$-function represents
(in the sense of the number of divisor points)
``a quarter" of the $\theta_1(u/\omega_2|\omega_1/\omega_2)$ Jacobi
theta function. Correspondingly, one can consider
equation \eqref{e-gamma-eq} as a composition of four equations for
 $\Gamma_3(u;\mathbf{\omega})$ with different arguments and quasiperiods
and represent the elliptic gamma functions as ratios of four Barnes gamma
functions of the third order with some simple exponential multipliers
\cite{FR,spi:theta2}. For some other important results for the generalized
gamma functions, see \cite{nar,rai:limits}.

\medskip {} {\bf \em The elliptic beta integral.}
It is convenient to use the compact notation
\begin{eqnarray*}
&& \eg(a_1,\ldots,a_k;p,q):=\eg(a_1;p,q)\cdots \eg(a_k;p,q),\quad
\\ &&
\Gamma(tz^{\pm 1};p,q):=\Gamma(tz;p,q)\Gamma(tz^{-1};p,q),\quad
\Gamma(z^{\pm2};p,q):=\Gamma(z^2;p,q)\Gamma(z^{-2};p,q)
\end{eqnarray*}
for working with elliptic hy\-per\-geo\-met\-ric integrals.
We start consideration from the elliptic beta integral discovered by the author
in \cite{spi:umn}.

\begin{theorem}
Take eight complex parameters $t_1,\ldots, t_6$,
and $p,q$, satisfying the constraints $|p|, |q|, |t_j|<1$ and $\prod_{j=1}^6 t_j=pq$.
Then the following equality is true
\be
\kappa\int_\T\frac{\prod_{j=1}^6
\eg(t_jz^{\pm 1};p,q)}{\eg(z^{\pm 2};p,q)}\frac{dz}{z}
=\prod_{1\leq j<k\leq6}\eg(t_jt_k;p,q),
\label{ell-int}\ee
where $\T$ denotes the positively oriented unit circle
and $\kappa=(p;p)_\infty(q;q)_\infty/4\pi\textup{i}.$
\end{theorem}
The first proof of this formula was based on the elliptic
extension of Askey's method \cite{ask:beta}.
A particularly short proof was given in \cite{spi:short}.
It is based on the partial $q$-difference equation
\begin{eqnarray*}
&& \rho(z;qt_1,t_2,\ldots,t_5;p,q)-\rho(z;t_1,\ldots,t_5;p,q)
\\ && \makebox[4em]{}
=g(q^{-1}z)\rho(q^{-1}z;,t_1,\ldots,t_5;p,q)-g(z)\rho(z;t_1,\ldots,t_5;p,q),
\nonumber\end{eqnarray*}
where $\rho(z;\underline{t};p,q)$ is the integral kernel divided by
the right-hand side expression in equality \eqref{ell-int} with $t_6$
replaced by $pq/t_1\cdots t_5$ and
$$
g(z)=
\frac{\prod_{k=1}^5\theta(t_kz;p)}{\prod_{k=2}^5\theta(t_1t_k;p)}
\frac{\theta(t_1\prod_{j=1}^5t_j;p)}
{\theta(z^2,z\prod_{j=1}^5t_j;p)}\frac{t_1}{z}.
$$
Dividing the above equation by $\rho(z;\underline{t};p,q)$, one comes to
a specific identity for elliptic functions.
A similar $p$-difference equation is obtained after permutation of $p$ and $q$.
Jointly they show that the integral $I(\underline{t})=\int_\T\rho(z;\underline{t};p,q)dz/z$
satisfies the equations $I(qt_1,t_2,\ldots,t_5)=I(pt_1,t_2,\ldots,t_5)=I(\underline{t})$.
In order to see this it is necessary to integrate the equations for
$\rho(z;\underline{t};p,q)$ over $z\in\mathbb{T}$ under the conditions
$|t_k|<1, k=1,\ldots,5$, and $\prod_{k=1}^5|t_k|>|p|, |q|$.
For incommensurate  $p$ and $q$ the invariance under scaling by these parameters
proves that the analytically continued function  $I(\underline{t})$
does not depend on $t_1$ and, in this way, is a constant independent on
all the parameters $t_j$.
Taking a special limit of parameters $t_j$ such that integral's value is asymptotically
given by the sum of residues of a fixed pair of poles (see below),
one finds this constant.

The elliptic beta integral \eqref{ell-int} defines the most general known
univariate exact integration formula generalizing Euler's beta integral.
For $p\to 0$, one obtains the Rahman integral \cite{rah:integral}
(see Theorem 10.8.2 in \cite{aar}), which reduces to the well known
Askey-Wilson $q$-beta integral \cite{aw} (see Theorem 10.8.1 in
\cite{aar}) if one of the parameters vanishes.
The binomial theorem $_1F_0(a;x)=(1-x)^{-a}$ (see formula (2.1.6) in \cite{aar})
was proved by Newton. The $q$-binomial theorem $_1\varphi_0(t;q;x)=(tx;q)_\infty/(x;q)_\infty$
(see Ch. 10.2 in \cite{aar}) was established by Gauss and several other mathematicians.
These formulas represent the simplest plain and $q$-hypergeometric function identities.
At the elliptic level, this role is played by the elliptic beta integral evaluation,
i.e. formula \eqref{ell-int} can be considered as an elliptic binomial theorem.

Replace in formula \eqref{ell-int}
$\T$ by a contour $C$ which separates sequences of the integrand poles
converging to zero along the points $z=t_jq^kp^m,\, k,m\in\Z_{\geq 0}$, from
their reciprocals obtained by the change $z\to 1/z$, which go to infinity.
This allows one to lift the constraints $|t_j|<1$ without changing the right-hand side
of formula \eqref{ell-int}. Substitute now $t_6=pq/A$, $A=\prod_{k=1}^5t_k$,
and suppose that $|t_m|<1,\, m=1,\ldots ,4,$ $|pt_5|<1<|t_5|$, $|pq|<|A|$,
and the arguments of $t_1,\ldots,t_5,$ and $p,q$ are linearly independent
over $\Z$. Then the following equality takes place \cite{die-spi:elliptic}:
\be
\kappa\int_C \Delta_E(z,\underline{t})\frac{d z}{z} =
\kappa\int_\mathbb{T} \Delta_E(z,\underline{t})\frac{d z}{z}
+c_0(\underline{t}) \sum_{|t_5q^n|>1,\, n \geq 0}
\nu_n(\underline{t}),
\lab{res}\ee
where $\Delta_E(z,\underline{t})=\prod_{m=1}^5\eg(t_mz^{\pm 1};p,q)
/\eg(z^{\pm2},Az^{\pm 1};p,q)$ and
\begin{eqnarray*} \nonumber
c_0(\underline{t}) =
\frac{\prod_{m=1}^4\eg(t_mt_5^{\pm 1};p,q)}
{\eg(t_5^{-2},A t_5^{\pm 1};p,q)}, \qquad
\nu_n(\underline{t}) =
 \frac{\theta(t_5^2q^{2n};p)}{\theta(t_5^2;p)}
\prod_{m=0}^5 \frac{\theta(t_mt_5)_n} {\theta(qt_m^{-1}t_5)_n}\, q^n.
\end{eqnarray*}
We have introduced here a new parameter $t_0$ with the help
of the relation $\prod_{m=0}^5t_m=q$ and used the elliptic
Pochhammer symbol
$$
\theta(t)_n=\prod_{j=0}^{n-1}\theta(tq^j;p)=\frac{\eg(tq^n;p,q)}{\eg(t;p,q)}, \qquad
\theta(t_1,\ldots,t_k)_n:=\prod_{j=1}^k\theta(t_j)_n
$$
(the indicated ratio of elliptic gamma functions defines $\theta(t)_n$
for arbitrary $n\in\C$). The multiplier $\kappa$ is absent
in the coefficient $c_0$ due to the relation
$\lim_{z\to 1}(1-z)\eg(z;p,q)=1/(p;p)_\infty(q;q)_\infty$
and doubling of the number of residues because of the symmetry $z\to z^{-1}$.

In the limit $t_5t_4\to q^{-N},\, N\in\Z_{\geq 0}$, the integral over the contour $C$
(equal to the right-hand side of equality \eqref{ell-int}) and the multiplier
$c_0(\underline{t})$ in front of the sum of residues
diverge, whereas the integral over the unit circle $\mathbb{T}$
remains finite. After dividing all the terms by $c_0(\underline{t})$
and going to the limiting relation, we obtain the Frenkel-Turaev
summation formula
\begin{equation}
\sum_{n=0}^N\nu_n(\underline{t})=
\frac{\theta(qt_5^2,\frac{q}{t_1t_2},
\frac{q}{t_1t_3},\frac{q}{t_2t_3})_N }
{\theta(\frac{q}{t_1t_2t_3t_5},\frac{qt_5}{t_1},
\frac{qt_5}{t_2},\frac{qt_5}{t_3})_N},
\label{ft-sum}\end{equation}
which was established for the first time in \cite{ft} by a completely
different method. For $N=0$ this equality trivializes and proves that
the integral considered earlier  $I(\underline{t})=1$.
For $p\to 0$ and fixed parameters, formula \eqref{ft-sum}
reduces to the Jackson sum for a terminating $_8\varphi_7$-series
(see Ex. 16 in Ch. 10 and formula (12.3.5) in \cite{aar}).
We stress that all terminating elliptic hypergeometric series
identities like identity \eqref{ft-sum} represent relations between ordinary
elliptic functions, i.e. they do not involve principally new special
functions in contrast to the elliptic hypergeometric integral identities.

\medskip {} {\bf \em General elliptic hy\-per\-geo\-met\-ric functions.}
Definitions of the general elliptic hy\-per\-geo\-met\-ric series and integrals
were  given and investigated in detail in \cite{spi:theta1}
and \cite{spi:theta2}, respectively.  So, a formal series $\sum_{n\in\Z}c_n$
is called an elliptic hy\-per\-geo\-met\-ric series if $c_{n+1}=h(n)c_n,$ where $h(n)$
is some elliptic function of $n\in\C$. This definition is
contained implicitly in the considerations of \cite{spi-zhe:spectral}.
It is well known \cite{akh}
that an arbitrary elliptic function $h(u)$ of order $s+1$ with the periods
$\omega_2/\omega_1$ and $\omega_3/\omega_1$ can be represented in the form
\begin{equation}
h(u)=y\;\prod_{k=1}^{s+1}\frac{\theta(t_kz;p)}{\theta(w_kz;p)},\qquad z=q^u.
\label{el-fn}\end{equation}
The equality $h(u+\omega_2/\omega_1)=h(u)$ is evident, and the periodicity
$h(u+\omega_3/\omega_1)=h(u)$ brings in the balancing condition
$
\prod_{k=1}^{s+1} t_k=\prod_{k=1}^{s+1} w_k.
$
Because of the factorization of $h(n)$, in order to determine the coefficients
$c_n$ it suffices to solve the equation
$a_{n+1}=\theta(tq^n;p)\, a_n,$ which leads to the elliptic Pochhammer symbol
$a_n=\theta(t)_n\, a_0$. The explicit form of the bilateral elliptic
hy\-per\-geo\-met\-ric series is now easily found to be
$$
{}_{s+1}G_{s+1}\bigg({t_1,\ldots,t_{s+1}\atop w_1,\ldots, w_{s+1}};q,p;y\bigg)
:=\sum_{n\in\Z}\prod_{k=1}^{s+1}\frac{\theta(t_k)_n}
{\theta(w_k)_n}\,y^n,
$$
where we have chosen the normalization $c_0=1$. By setting
$w_{s+1}=q, \, t_{s+1}=: t_0$, we obtain the one sided series
\begin{equation}
{}_{s+1}E_s\bigg({t_0,t_1,\ldots  ,t_s\atop w_1,\ldots  ,w_s};q,p;y\bigg)
:=\sum_{n\in\Z_{\geq 0}}\frac{\theta(t_0,t_1,\ldots, t_s)_n}{
\theta(q,w_1,\ldots,w_s)_n}\, y^n.
\label{E-series}\end{equation}
For fixed $t_j$ and $w_j$, the function
${}_{s+1}E_s$ degenerates in the limit $p\to 0$ to the basic $q$-hy\-per\-geo\-met\-ric series
${}_{s+1}\varphi_s$ satisfying the condition $\prod_{k=0}^st_s=q\prod_{k=1}^sw_s$.
The infinite series \eqref{E-series} does not converge in general, and we therefore assume
its termination due to the condition $t_k=q^{-N}p^M$ for some $k$ and
$N\in\Z_{\geq 0},\, M\in\Z.$ The additive system of notation for these series
(see, e.g., Ch. 11 in \cite{gas-rah:basic} or \cite{spi:thesis})
is more convenient for consideration of certain questions, but we skip it here.

The series \eqref{E-series} is called well-poised if
$t_0q=t_1w_1=\ldots=t_sw_s$. In this case the balancing condition
takes the form $t_1\cdots t_s=\pm q^{(s+1)/2}t_0^{(s-1)/2}$, and
the functions $h(u)$ and $_{s+1}E_s$ become invariant under the changes
$t_j\to pt_j,\, j=1,\ldots,s-1,$ and $t_0\to p^2t_0$. For odd $s$ and
balancing condition of the form $t_1\cdots t_s=+q^{(s+1)/2}t_0^{(s-1)/2}$,
one has the symmetry $t_0\to pt_0$ and $_{s+1}E_s$ becomes an elliptic function
of all free parameters $\log t_j,\, j=0,\ldots, s-1,$ with equal periods
(such functions were called in \cite{spi:theta1,spi:thesis} totally
elliptic functions). Under the four additional constraints
$t_{s-3}=q\sqrt{t_0},\, t_{s-2}=-q\sqrt{t_0},\, t_{s-1}=
q\sqrt{t_0/p},\, t_s=-q\sqrt{pt_0}$, connected to doubling of
the argument of theta functions, the series are called very-well-poised.
In \cite{spi:bailey1}, a special notation was introduced
for the very-well-poised elliptic hy\-per\-geo\-met\-ric series:
\begin{eqnarray}\label{V-series}
&&
{}_{s+1}E_{s}\bigg({t_0,t_1,\ldots,t_{s-4},q\sqrt{t_0},-q\sqrt{t_0},
q\sqrt{t_0/p},-q\sqrt{pt_0}  \atop
qt_0/t_1,\ldots  ,qt_0/t_{s-4},\sqrt{t_0}, -\sqrt{t_0},\sqrt{pt_0},
-\sqrt{t_0/p}};q,p;-y\bigg)
\\ && \makebox[2em]{}
= \sum_{n=0}^\infty
\frac{\theta(t_0q^{2n};p)}{\theta(t_0;p)}\prod_{m=0}^{s-4}
\frac{\theta(t_m)_n}{\theta(qt_0t_m^{-1})_n}(qy)^n =:
{}_{s+1}V_{s}(t_0;t_1,\ldots,t_{s-4};q,p;y),
\nonumber\end{eqnarray}
where the balancing condition has the form
$\prod_{k=1}^{s-4}t_k=\pm t_0^{(s-5)/2}q^{(s-7)/2}$,
and for odd $s$ we assume the positive sign choice for preserving
the symmetry $t_0\to pt_0$. If $y=1$, then $y$
is omitted in the series notation. Summation formula
\eqref{ft-sum} gives thus a closed form expression for the terminating
$_{10}V_9(t_0;t_1,\ldots,t_5;q,p)$-series.

A contour integral $\int_C\Delta(u)du$ is called an elliptic hy\-per\-geo\-met\-ric
integral if its kernel $\Delta(u)$ satisfies the system of three equations
\be
\Delta(u+\omega_k)=h_k(u)\Delta(u),\quad k=1,2,3,
\lab{ell-gip}\ee
where $\omega_{1,2,3}\in \C$ are some pairwise incommensurate parameters
and $h_k(u)$ are some elliptic functions with periods $\omega_{k+1},$ $\omega_{k+2}$
(we set $\omega_{k+3}=\omega_k$). One can weaken the requirement \eqref{ell-gip}
by keeping only one equation, but then there appears a functional freedom
in the choice of $\Delta(u)$, which should be fixed in some other way.

Omitting the details of such considerations from \cite{spi:theta2,spi:thesis},
we present the general form of permissible functions $\Delta(u)$. We suppose that this
function satisfies the equations \eqref{ell-gip} for $k=1,2$, where
$$
h_1(u)=y_1\prod_{j=1}^s\frac{\theta(t_je^{2\pi\textup{i}u/\omega_2};p)}
{\theta(w_je^{2\pi\textup{i}u/\omega_2};p)}, \quad
h_2(u)=y_2\prod_{j=1}^{\ell }\frac{\theta(\tilde t_je^{-2\pi\textup{i}u/\omega_1};r)}
{\theta(\tilde w_je^{-2\pi\textup{i}u/\omega_1};r)},
$$
$|p|,|r|<1$ and $\prod_{j=1}^st_j=\prod_{j=1}^s w_j$,
$\prod_{j=1}^{\ell }\tilde t_j=\prod_{j=1}^{\ell }\tilde w_j$.
If we take $|q|<1$, then the most general meromorphic $\Delta(u)$
has the form
\begin{equation}\label{delta-ell}
\Delta(u)=\prod_{j=1}^s\frac{\eg(t_je^{2\pi\textup{i}\frac{u}{\omega_2}};p,q)}
{\eg(w_je^{2\pi\textup{i}\frac{u}{\omega_2}};p,q)}
\prod_{j=1}^{\ell }\frac{\Gamma(\tilde t_je^{-2\pi\textup{i}\frac{u}{\omega_1}};
\tilde q,r)} {\Gamma(\tilde w_je^{-2\pi\textup{i}\frac{u}{\omega_1}}; \tilde q,r)}
\prod_{k=1}^m \frac{\theta(a_ke^{2\pi\textup{i}\frac{u}{\omega_2}};q)}
{\theta(b_ke^{2\pi\textup{i}\frac{u}{\omega_2}};q)}\, e^{cu+d},
\end{equation}
where the parameters $d\in\C$ and $m\in\Z_{\geq 0}$ are arbitrary, and $a_k, b_k, c$ are connected with
$y_1$ and $y_2$ by the relations $y_2=e^{c\omega_2}$ and
$y_1=e^{c\omega_1} \prod_{k=1}^mb_ka_k^{-1}$. It appears that the function
$h_3(u)$ cannot be arbitrary -- it is determined from the integral kernel \eqref{delta-ell}.

For $|q|=1$ it is necessary to choose $\ell =s$ in formula \eqref{delta-ell}
and fix parameters in such a way that the $\Gamma$-functions are combined
to the modified elliptic gamma function $G(u;\mathbf{\omega})$
(it is precisely in this way that this function was built in \cite{spi:theta2}):
\begin{equation}\label{delta-ell-unit}
\Delta(u)=\prod_{j=1}^s\frac{G(u+g_j;\mathbf{\omega})}
{G(u+v_j;\mathbf{\omega})}\, e^{cu+d},
\end{equation}
where the parameters $g_j,\, v_j$ are connected to $t_j,\, w_j$ by the
relations $t_j=e^{2\pi\textup{i}g_j/\omega_2},$  $w_j=e^{2\pi\textup{i} v_j/\omega_2}$,
and $y_{1,2}=e^{c\omega_{1,2}}$. The integrals $\int_C\Delta(u) du$
with kernels of the indicated form define elliptic analogues of the
Meijer function. For even more general theta hy\-per\-geo\-met\-ric
integrals, see \cite{spi:theta2}.

We limit consideration to the case when both $\ell$ and $m$ in
\eqref{delta-ell} are equal to zero.
The corresponding integrals are called well-poised, if
$t_1w_1=\ldots=t_sw_s=pq$. The additional condition of very-well-poisedness
fixes eight parameters
$t_{s-7},\ldots,$ $t_s=\{\pm(pq)^{1/2},$ $\pm q^{1/2}p,$ $\pm p^{1/2}q,\pm pq\}$
and doubles the argument of the elliptic gamma function:
$\prod_{j=s-7}^{s}\eg(t_jz;p,q)=1/\eg(z^{-2};p,q).$
The most interesting are the very-well-poised elliptic hy\-per\-geo\-met\-ric
integrals with even number of parameters
\begin{equation}
I^{(m)}(t_1,\ldots,t_{2m+6})
=\kappa\int_\T\frac{\prod_{j=1}^{2m+6}\eg(t_jz^{\pm 1};p,q)}{\eg(z^{\pm2};p,q)}
\frac{dz}{z}, \quad \prod_{j=1}^{2m+6}t_j=(pq)^{m+1},
\label{gen-int}\end{equation}
with $|t_j|<1$ and ``correct" choice of the sign in the
balancing condition. They represent integral analogues of the
$_{s+1}V_s$-series with odd $s$, ``correct" balancing condition and
the argument $y=1$, in the sense that such series appear as residue
sums of particular pole sequences of the kernel of $I^{(m)}$.
Note that $I^{(0)}$ coincides with the elliptic beta integral.

Properties of the elliptic functions explain the origins of
hy\-per\-geo\-met\-ric notions of balancing, well-poisedness,
and very-well-poisedness.
However, strictly speaking these notions are consistently defined
only at the elliptic level, because there are limits to such $q$-hy\-per\-geo\-met\-ric
identities in which they are not preserved any more
\cite{rai:abelian,spi:bailey1,war:summation2}! The fact of unique determination of the
balancing condition for series \eqref{V-series} with odd $s$ and integrals
\eqref{gen-int} (precisely these objects emerge in interesting applications)
illustrates a deep internal tie between the ``elliptic" and ``hy\-per\-geo\-met\-ric"
classes of special functions. Multivariable elliptic hy\-per\-geo\-met\-ric series
and integrals are defined analogously to the univariate case -- it
is necessary to use systems of finite difference equations for
kernels with the coefficients given by elliptic functions of all
summation or integration variables \cite{spi:theta1,spi:theta2},
which is a natural generalization of the approach of Pochhammer and Horn
to functions of hy\-per\-geo\-met\-ric type \cite{aar,ggr}.

\medskip {} {\bf \em An elliptic analogue of the Euler-Gauss hy\-per\-geo\-met\-ric function.}
Take eight parameters $t_1,\ldots, t_8\in\C $ and two base variables
$p,q\in\C$ satisfying the constraints $ |p|,|q|<1$ and
$\prod_{j=1}^8 t_j=p^2q^2$ (the balancing condition).
For all $|t_j|<1$ an elliptic analogue of the Euler-Gauss
hy\-per\-geo\-met\-ric function $_2F_1(a,b;c;x)$ (see Ch. 2 in \cite{aar})
 is defined by the integral \cite{spi:thesis}
\be
V(\underline{t})\equiv V(t_1,\ldots,t_8;p,q):=\kappa\int_\T\frac{\prod_{j=1}^8
\eg(t_jz^{\pm 1};p,q)}{\eg(z^{\pm 2};p,q)}\frac{dz}{z},
\lab{ehf}\ee
i.e. by the choice $m=1$ in expression \eqref{gen-int}.
Note that it can be reduced to both Euler and Barnes type integral
representations of $_2F_1$-series.
For other admissible values of parameters, the $V$-function is defined
by the analytical continuation of expression \eqref{ehf}.
From this continuation one can see that the $V$-function is meromorphic
for all values of parameters $t_j\in\C^*$ when the contour of integration
is not pinched. To see this, compute residues of the integrand
poles and define the analytically continued function as a sum of the integral over some fixed
contour and residues of the poles crossing this contour.
More precisely, $\prod_{1\leq j<k\leq 8}(t_jt_k;p,q)_\infty V(\underline{t})$
is a holomorphic function of parameters \cite{rai:trans}. As shown in \cite{S5},
the $V$-function has delta-function type singularities at certain values of $t_j$'s.

The first nontrivial property of function \eqref{ehf} consists in its reduction to
the elliptic beta integral under the condition for a pair of parameters
$t_jt_k=pq$, $j\neq k$ (expression  \eqref{ell-int} appears from $t_7t_8=pq$).
The $V$-function is evidently symmetric in $p$ and $q$. It is invariant
also under the $S_8$-group of permutations of parameters $t_j$
isomorphic to the Weyl group $A_7$.
Consider the double integral
$$
\kappa \int_{\T^2}
\frac{\prod_{j=1}^4\eg(a_jz^{\pm 1}, b_jw^{\pm 1};p,q)\;\eg(cz^{\pm 1} w^{\pm 1};p,q)}
{\eg(z^{\pm 2},w^{\pm 2};p,q)} \frac{dz}{z}\frac{dw}{w},
$$
where $a_j,b_j,c\in \C $, $|a_j|,|b_j|,|c|<1,$ and
$c^2\prod_{j=1}^4a_j= c^2\prod_{j=1}^4b_j=pq.$
Using formula \eqref{ell-int} for integration over $z$ or $w$
(the permutation of the order of integrations is permitted), we obtain
the following transformation formula:
\be
V(\underline{t})=\prod_{1\le j<k\le 4}\eg(t_jt_k,t_{j+4}t_{k+4};p,q)\,
V(\underline{s}),
\lab{E7-1}\ee
where $|t_j|, |s_j|<1$, and
$$
\left\{
\begin{array}{cl}
s_j =\rho^{-1} t_j,&   j=1,2,3,4  \\
s_j = \rho t_j, &    j=5,6,7,8
\end{array}
\right.;
\quad \rho=\sqrt{\frac{t_1t_2t_3t_4}{pq}}=\sqrt{\frac{pq}{t_5t_6t_7t_8}}.
$$
This fundamental relation was derived
by the author in \cite{spi:theta2}, where the function $V(\underline{t})$
appeared for the first time. It represents an elliptic analogue
(moreover, integral generalization) of Bailey's transformation
for four non-terminating  $_{10}\varphi_9$-series \cite{gas-rah:basic}.

Repeat transformation \eqref{E7-1} once more with the parameters $s_{3,4,5,6}$,
playing the role of $t_{1,2,3,4}$, and permute parameters $t_3,t_4$ with $t_5,t_6$
in the resulting expression. This yields the relation
\be
V(\underline{t})=\prod_{j,k=1}^4
\eg(t_jt_{k+4};p,q)\ V(T^{1\over 2}\!/t_1,\ldots,T^{1\over 2}\!/t_4,
U^{1\over 2}\!/t_5,\ldots,U^{1\over 2}\!/t_8),
\lab{E7-2}\ee
where $ T=t_1t_2t_3t_4$, $ U=t_5t_6t_7t_8$ and
$|T|^{1/2}<|t_j|<1,$ $|U|^{1/2}<|t_{j+4}|<1,\, j=1,2,3,4$.
Now equating the right-hand sides of relations \eqref{E7-1} and \eqref{E7-2},
and expressing parameters $t_j$ in terms of $s_j$, one obtains the third relation
\be
V(\underline{s})=\prod_{1\le j<k\le 8}\eg(s_js_k;p,q)\,
V(\sqrt{pq}/s_1,\ldots,\sqrt{pq}/s_8),
\lab{E7-3}\ee
where $|pq|^{1/2}<|s_j|<1$ for all  $j$.

Consider the Euclidean space $\R^8$ with the scalar product
$\langle x, y\rangle$ and an orthonormal basis $e_i\in \R^8$,
$\langle e_i, e_j\rangle=\delta_{ij}.$ The root system $A_7$ consists of
the vectors $v=\{e_i-e_j,\, i\neq j\}.$ Its Weyl group consists of
the reflections $x\to S_v(x)=x-2v\langle v, x\rangle/\langle v, v\rangle$
acting in the hyperplane orthogonal to the vector $\sum_{i=1}^8e_i$
(i.e., the coordinates of the vectors $x=\sum_{i=1}^8x_ie_i$ satisfy
the constraint $\sum_{i=1}^8x_i=0$), and it coincides with the
permutation group $S_8$.

Connect parameters of the $V(\underline{t})$-function to the
coordinates $x_j$ as $t_j=e^{2\pi\textup{i} x_j}(pq)^{1/4}$.
Then the balancing condition assumes the form $\sum_{i=1}^8 x_i=0$.
The first $V$-function transformation \eqref{E7-1} is now easily seen to
correspond to the reflection $S_v(x)$ for the vector
$v=(\sum_{i=5}^8e_i-\sum_{i=1}^4 e_i)/2$ having the canonical length
$\langle v, v\rangle=2$. This reflection extends the group $A_7$ to the
exceptional Weyl group $E_7$. Relations \eqref{E7-2} and \eqref{E7-3} were
proved in a different fashion by Rains in \cite{rai:trans}, where it was
indicated that these transformations belong to the group $E_7$.

Denote by $V(qt_j,q^{-1}t_k)$ elliptic hy\-per\-geo\-met\-ric functions
 contiguous to $V(\underline{t})$ in the sense that
 $t_j$ and $t_k$ are replaced by $qt_j$ and $q^{-1}t_k$, respectively.
The following contiguous relation for the $V$-functions is valid
\be
t_7\theta\left(t_8t_7^{\pm 1}/q;p\right)V(qt_6,q^{-1}t_8)
-(t_6\leftrightarrow t_7)=t_7\theta\left(t_6t_7^{\pm 1};p\right) V(\underline{t}),
\lab{con_1}\ee
where $(t_6\leftrightarrow t_7)$ denotes the permutation of parameters in
the preceding expression (such a relation was used already in \cite{spi:umn}).
Indeed, for $y=t_6, w=t_7, $ and $x=q^{-1}t_8$ the Riemann relation
\eqref{ident} is equivalent to the $q$-difference equation for $V$-function's
integrand
 $\Delta(z,\underline{t})=\prod_{k=1}^8\eg(t_kz^{\pm 1};p,q)/\eg(z^{\pm2};p,q)$
coinciding with \eqref{con_1} after replacement of $V$-functions
by $\Delta(z,\underline{t})$ with appropriate parameters.
Integration of this equation over the contour $\T$ yields formula \eqref{con_1}.
Substitute now the symmetry transformation \eqref{E7-3} in  \eqref{con_1}
and obtain the second contiguous relation
\begin{eqnarray*}
t_6\theta\Bigl(\frac{t_7}{qt_8};p\Bigr)\prod_{k=1}^5
\theta\Bigl(\frac{t_6t_k}{q};p\Bigr)V(q^{-1}t_6,qt_8)
-(t_6\leftrightarrow t_7)
=t_6\theta\Bigl(\frac{t_7}{t_6};p\Bigr)
\prod_{k=1}^5\theta(t_8t_k;p) V(\underline{t}).
\end{eqnarray*}
An appropriate combination of these two equalities yields the equation
\be
\mathcal{A}(\underline{t})\Big(U(qt_6,q^{-1}t_7)-U(\underline{t})\Big)
+(t_6\leftrightarrow t_7)+ U(\underline{t})=0,
\lab{pre-ehf}\ee
where we have denoted
$U(\underline{t})=V(\underline{t})/\eg(t_6t_8^{\pm 1},t_7t_8^{\pm 1};p,q)$ and
\be
\mathcal{A}(\underline{t})=\frac{\theta(t_6/qt_8,t_6t_8,t_8/t_6;p)}
                 {\theta(t_6/t_7,t_7/qt_6,t_6t_7/q;p)}
\prod_{k=1}^5\frac{\theta(t_7t_k/q;p)}{\theta(t_8t_k;p)}.
\ee
Substituting $t_j=e^{2\pi\textup{i} g_j/\omega_2}$, one can check that the potential
$\mathcal{A}(\underline{t})$ is a modular invariant elliptic function of
the variables $g_1,\ldots,g_7$, i.e. it does not change after the replacements
$g_j\to g_j+\omega_{2,3}$ or $(\omega_2,\omega_3)\to(-\omega_3,\omega_2)$.

Now denote $t_6=cx,\ t_7=c/x,$ and introduce new variables
$$
\ve_k=\frac{q}{ct_k},\; k=1,\ldots,5,\quad \ve_8=\frac{c}{t_8},
\quad \ve_7=\frac{\ve_8}{q},\quad c=\frac{\sqrt{\ve_6\ve_8}}{p^2}.
$$
In terms of $\ve_k$ the balancing condition
takes the standard form $\prod_{k=1}^8\ve_k=p^2q^2$.
After the replacement of $U(\underline{t})$ in formula \eqref{pre-ehf}
by some unknown function $f(x)$, we obtain a $q$-difference
equation of the second order which is called
{\em the elliptic hy\-per\-geo\-met\-ric equation} \cite{spi:thesis,spi:cs}:
\begin{eqnarray}
&& A(x)\left( f(qx)-f(x)\right)
+ A(x^{-1})\left( f(q^{-1}x)-f(x)\right) + \nu f(x)=0,
\label{eheq}
\\ && \qquad
A(x)=\frac{\prod_{k=1}^8 \theta(\ve_kx;p)}{\theta(x^2,qx^2;p)},
\qquad
\nu=\prod_{k=1}^6\theta\left(\frac{\ve_k\ve_8}{q};p\right).
\lab{pot}\end{eqnarray}
We have already one functional solution of this equation
\be
f_1(x)=\frac{ V(q/c\ve_1,\ldots,q/c\ve_5,cx,c/x,c/\ve_8;p,q)}
{\eg(c^2x^{\pm 1}/\ve_8,x^{\pm 1} \ve_8;p,q)},
\lab{sol1}\ee
where it is necessary to impose the constraints
(in the previous parametrization)
$\sqrt{|pq|}<|t_j|<1,\, j=1,\ldots, 5,$ and $\sqrt{|pq|}<|q^{\pm1}t_6|,|q^{\pm1}t_7|,
|q^{\pm1}t_8|<1$, which can be relaxed by analytical continuation.
Other independent solutions can be obtained by the multiplication of one
of the parameters $\ve_1,\ldots,\ve_5$, and $x$ by powers of $p$ or
by permutations of $\ve_1,\ldots,\ve_5$ with $\ve_6$.

Denote $\ve_k=e^{2\pi\textup{i} a_k/\omega_2}$,
$x=e^{2\pi\textup{i} u/\omega_2}$, and
$F_1(u;\underline{a};\omega_1,\omega_2,\omega_3):=f_1(x)$.
Then one can check that equation
\eqref{eheq} is invariant with respect to the modular
transformation $(\omega_2,\omega_3)\to (-\omega_3,\omega_2)$.
Therefore one of the linear independent solutions of \eqref{eheq}
has the form $F_2(u;\underline{a};\omega_1,\omega_2,\omega_3):=
F_1(u;\underline{a};\omega_1,-\omega_3,\omega_2).$
The same solution would be obtained if we repeat the derivation of
equation  \eqref{eheq} and its solution \eqref{sol1}
after replacing $\Gamma$-functions by the modified elliptic
gamma function $G(u;\mathbf{\omega})$. This shows that $F_2$-function
is well defined even for $|q|=1$. Different limiting transitions from
the $V$-function and other elliptic hy\-per\-geo\-met\-ric
integrals to $q$-hy\-per\-geo\-met\-ric integrals of the Mellin-Barnes or
Euler type are described in \cite{spi:thesis,spi:cs} and much more
systematically in \cite{BultPhD,br,brs,rai:limits}.

\medskip {} {\bf \em Biorthogonal functions of the hy\-per\-geo\-met\-ric type.}
In analogy with the residue calculus for the elliptic beta integral
 \eqref{res}, one can consider the sum of residues for a particular
geometric progression of poles of the $V$-function kernel for one of the
parameters. This leads to the very-well-poised ${}_{12}V_{11}$-elliptic
hy\-per\-geo\-met\-ric series the termination of which
is guaranteed by a special discretization of the chosen parameter.
In this way one can rederive contiguous relations for the terminating
${}_{12}V_{11}$-series of \cite{spi-zhe:spectral,spi-zhe:rims}
out of the contiguous relations for the $V$-function, which
we omit here. For instance, this yields the following particular
solution of the elliptic hy\-per\-geo\-met\-ric equation \eqref{eheq}:
\be
R_n(x;q,p)=
{}_{12}V_{11}\left(\frac{\ve_6}{\ve_8};\frac{q}{\ve_1\ve_8},
\frac{q}{\ve_2\ve_8},\frac{q}{\ve_3\ve_8},
\frac{qp}{\ve_4\ve_8},\frac{qp}{\ve_5\ve_8},\ve_6x,\frac{\ve_6}{x};q,p\right),
\lab{R_n}\ee
where $pq/\ve_4\ve_8=q^{-n},$ $n\in\Z_{\geq 0}$
(we recall that $\prod_{k=1}^8\ve_k=p^2q^2$).
Properties of the $R_n$-function were described in  \cite{spi:theta2},
whose notation passes to ours after the replacements
$t_{0,1,2}\to \ve_{1,2,3},\, t_3\to\ve_6,\,t_4\to\ve_8,\,
\mu\to\ve_4\ve_8/pq$, and $A\mu/qt_4\to pq/\ve_5\ve_8$.

Equation \eqref{eheq} is symmetric in $\ve_1,\ldots,\ve_6$.
The series \eqref{V-series} is elliptic in all parameters, therefore function
\eqref{R_n} is symmetric in $\ve_1,\ldots,\ve_5$ and each of these variables can be used
for terminating the series. A permutation  of $\ve_1,\ve_2,\ve_3,\ve_5$ with
$\ve_6$ yields $R_n(z;q,p)$ up to some multiplier independent
on $x$ due to an elliptic analogue of the Bailey transformation
for terminating ${}_{12}V_{11}$-series \cite{ft}, which
can be obtained by degeneration from equality \eqref{E7-1}.

The same contiguous relations for the $_{12}V_{11}$-series yield the
following three term recurrence relation for $R_n(x;q,p)$ in the index $n$:
\begin{eqnarray} \lab{ttr}
&& \makebox[-2em]{}
(z(x)-\alpha_{n+1})\rho(Aq^{n-1}/\ve_8)\left(R_{n+1}(x;q,p)-R_n(x;q,p)\right)
+(z(x)-\beta_{n-1}) \\ &&
\times \rho(q^{-n})\left(R_{n-1}(x;q,p)-R_n(x;q,p)\right)
+\delta (z(x)-z(\ve_6)) R_n(x;q,p)=0,
\nonumber \end{eqnarray}
where
\begin{eqnarray} &&
z(x)= \frac{\theta(x\xi^{\pm 1};p)}{\theta(x\eta^{\pm 1};p)}, \qquad
\alpha_n= z(q^n/\ve_8),\qquad \beta_n= z(Aq^{n-1}),
\nonumber  \\ &&
\rho(t)=\frac{\theta\left(t,\frac{\ve_6}{\ve_8t},\frac{q\ve_6}{\ve_8t},
\frac{qt}{\ve_1\ve_2},\frac{qt}{\ve_2\ve_3},\frac{qt}{\ve_1\ve_3},
\frac{q^2t\eta^{\pm 1}}{A};p\right)}
{\theta\left(\frac{qt^2\ve_8}{A},\frac{q^2t^2\ve_8}{A};p\right)},
\nonumber  \\ &&
\delta=\theta\left(\frac{q^2\ve_6}{A},\frac{q}{\ve_1\ve_8},
\frac{q}{\ve_2\ve_8},\frac{q}{\ve_3\ve_8},\ve_6\eta^{\pm 1};p\right).
\nonumber\end{eqnarray}
Here $A=\ve_1\ve_2\ve_3\ve_6\ve_8$, and $\xi$ and $\eta$ are arbitrary
gauge parameters, $\xi\neq \eta^{\pm1} p^k,\, k\in\Z$. The initial conditions
 $R_{-1}=0$ and $R_0=1$ guarantee that all the dependence on the variable
$x$ enters only through $z(x)$, and that $R_n(x)$ is a rational function of
$z(x)$ with poles at the points $\alpha_1,\ldots,\alpha_n$.

The elliptic hy\-per\-geo\-met\-ric equation for the $R_n$-function can be rewritten
in the form of a generalized eigenvalue problem
$\mathcal{D}_1R_n=\lambda_n\mathcal{D}_2R_n$ for some $q$-difference
operators of the second order $\mathcal{D}_{1,2}$ and discrete
spectrum  $\lambda_n$ \cite{spi:theta2}. We denote by $\phi_\lambda$ solutions
of an abstract spectral problem
$\mathcal{D}_1\phi_\lambda= \lambda \mathcal{D}_2\phi_\lambda$,
and by $\psi_\lambda$ solutions of the equation $\mathcal{D}_1^T\psi_\lambda=
\lambda \mathcal{D}_2^T\psi_\lambda$, where $\mathcal{D}_{1,2}^T$
are the operators conjugated with respect to some inner product
$\langle \psi|\phi\rangle$, i.e. $\langle \mathcal{D}_{1,2}^T\psi|\phi\rangle
=\langle \psi|\mathcal{D}_{1,2}\phi\rangle.$ Then
$0=\langle \psi_\mu|(\mathcal{D}_1-\lambda \mathcal{D}_2)\phi_\lambda\rangle
=(\mu-\lambda)\langle \mathcal{D}_2^T\psi_\mu|\phi_\lambda\rangle$,
i.e. the function $\mathcal{D}_2^T\psi_\mu$ is orthogonal to $\phi_\lambda$
for $\mu\neq\lambda$. As shown by Zhedanov \cite{zhe:gevp} (see also
\cite{spi-zhe:spectral,spi-zhe:rims}), this simple fact can be used for a
formulation of the theory of biorthogonal rational functions generalizing
orthogonal polynomials.
Analogues of the functions $\mathcal{D}_2^T\psi_\mu$ for $R_n(z;q,p)$ have the form
\be
T_n(x;q,p)=
{}_{12}V_{11}\left(\frac{A\ve_6}{q};\frac{A}{\ve_1},\frac{A}{\ve_2},\frac{A}{\ve_3},
\ve_6x,\frac{\ve_6}{x},\frac{qp}{\ve_4\ve_8},\frac{qp}{\ve_5\ve_8};q,p\right),
\lab{T_n}\ee
which are rational functions of $z(x)$ with poles at the points
$\beta_1,\ldots,\beta_n$.

Denote $R_{nm}(x):= R_n(x;q,p)R_m(x;p,q)$ and
$T_{nm}(x):= T_n(x;q,p)T_m(x;p,q)$, where all the $_{12}V_{11}$-series
terminate simultaneously because of the modified termination condition
$\ve_4\ve_8=p^{m+1}q^{n+1},\, n,m\in\Z_{\geq 0}$. The functions $R_{nm}$ now
solve not one but two generalized eigenvalue problems which differ from each
other by the permutation of  $p$ and $q$.

\begin{theorem}
The following two-index biorthogonality relation is true:
\be
\kappa\int_{C_{mn,kl}}T_{nl}(x)R_{mk}(x)
\frac{\prod_{j\in S}\eg(\ve_j x^{\pm 1};p,q)}
{\eg(x^{\pm2}, Ax^{\pm 1};p,q)}\frac{dx}{x}
=h_{nl}\: \delta_{mn}\: \delta_{kl},
\lab{2ib}\ee
where $S=\{1,2,3,6,8\}$, $C_{mn,kl}$ denotes the contour separating
sequences of points
$
x=\ve_jp^aq^b \,(j=1,2,3,6),\; \ve_8 p^{a-k}q^{b-m}, p^{a+1-l}q^{b+1-n}/A,\;
a,b\in \Z_{\geq 0},
$
from their $x\to x^{-1}$ reciprocals, and the normalization constants
have the form
\begin{eqnarray*}
h_{nl}&=&
\frac{\prod_{j< k,\, j,k\in S} \eg(\ve_j\ve_k;p,q)}
{\prod_{j\in S} \eg(A\ve_j^{-1};p,q)}\,h_n(q,p)\cdot h_l(p,q), \\
h_n(q,p)&=&\frac{\theta(A/q\ve_8;p)
\theta(q,q\ve_6/\ve_8,\ve_1\ve_2,\ve_1\ve_3,\ve_2\ve_3,A\ve_6)_n\,q^{-n}}
{\theta(Aq^{2n}/q\ve_8;p) \theta(1/\ve_6\ve_8,\ve_1\ve_6,\ve_2\ve_6,\ve_3\ve_6,
A/q\ve_6,A/q\ve_8)_n}.
\end{eqnarray*}
\end{theorem}

This theorem was proved in \cite{spi:theta2} by direct computation of
the integral in the left-hand side with the help of formula \eqref{ell-int}.
The appearance of the two-index orthogonality relations for functions of one variable
is a new phenomenon in the theory of special functions. It should be remarked
that only for $k=l=0$ there exists the limit $p\to 0$ and the resulting functions
$R_n(x;q,0)$, $T_n(x;q,0)$ coincide with
Rahman's family of continuous ${}_{10}\varphi_9$-biorthogonal rational functions
\cite{rah:integral}. A special limit $\text{Im}(\omega_3)\to\infty$
in the modular transformed $R_{nm}$ and $T_{nm}$ leads to the two-index
biorthogonal functions which are expressed as products of two modular conjugated
$_{10}\varphi_9$-series \cite{spi:thesis}. A special restriction for one of the
parameters in $R_n(x;q,p)$ and $T_n(x;q,p)$ leads to the biorthogonal rational
functions of a discrete argument derived by Zhedanov and the author in
\cite{spi-zhe:spectral} which generalizes
Wilson's functions \cite{wil:orthogonal}. All these functions are natural
generalizations of the Askey-Wilson polynomials \cite{aw}.

Note that $R_{nm}(x)$ and $T_{nm}(x)$ are meromorphic functions of the variable
$x\in\C^*$ with essential singularities at $x=0, \infty$ and only for
$k=l=0$ or $n=m=0$ do they become rational functions of some argument
depending on $x$. The continuous parameters biorthogonality relation
for the $V$-function itself was established in \cite{S5}. The
biorthogonal functions generated by the three-term recurrence
relation \eqref{ttr} after shifting $n$ by an arbitrary (complex)
number are not investigated yet. A generalization  of the described
``classical" biorthogonal functions to the ``semiclassical" level
associated with the higher order elliptic beta integrals \eqref{gen-int}
was suggested by Rains in \cite{rai:painleve}.

\medskip {} {\bf \em Elliptic beta integrals on root systems.}
Define a $C_n$ (or $BC_n$) root system analogue of
the constant $\kappa$:
$\kappa_{n}=(p;p)_\infty^n(q;q)_\infty^n/(2\pi\textup{i})^n 2^nn!$.
Describe now a $C_n$-elliptic beta integral representing a multiparameter
generalization of integral \eqref{ell-int}, which was classified in
\cite{die-spi:selberg} as an integral of type I.
\begin{theorem}\label{C-I}
Take $n$ variables $z_1,\ldots,z_n\in\T$ and complex parameters
$t_1,\ldots,$ $t_{2n+4}$ and $p,q$ satisfying the constraints
$|p|, |q|, |t_j|<1$ and $\prod_{j=1}^{2n+4}t_j=pq$. Then
\begin{eqnarray}\nonumber
&& \kappa_n\int_{\T^n}\prod_{1\leq j<k\leq n}\frac{1}{\eg(z_j^{\pm 1} z_k^{\pm 1};p,q)}
\prod_{j=1}^n\frac{\prod_{m=1}^{2n+4}\eg(t_mz_j^{\pm 1};p,q)}
{\eg(z_j^{\pm2};p,q)}\frac{dz_1}{z_1}\cdots\frac{dz_n}{z_n}
\\ && \makebox[8em]{}
=\prod_{1\leq m<s\leq 2n+4}\eg(t_mt_s;p,q).
\label{C-typeI}\end{eqnarray}
\end{theorem}
Formula \eqref{C-typeI} was suggested and partially confirmed
by van Diejen and the author in  \cite{die-spi:selberg}. It was
proved by different methods in \cite{rai:trans,RS,spi:thesis,spi:short}.
It reduces to one of Gustafson's integration
formulas \cite{gus:some} in a special $p\to 0$ limit.

\begin{theorem}
Take complex parameters $t, t_1,\ldots, t_6, p$ and $q$ restricted
by the conditions $|p|, |q|,$ $|t|,$ $|t_m| <1$ and $t^{2n-2}\prod_{m=1}^6t_m=pq$.
Then,
\begin{eqnarray}\nonumber
\kappa_n\int_{\T^n} \prod_{1\leq j<k\leq n}
\frac{\eg(tz_j^{\pm 1} z_k^{\pm 1};p,q)}{\eg(z_j^{\pm 1} z_k^{\pm 1};p,q)}
\prod_{j=1}^n\frac{\prod_{m=1}^6\eg(t_mz_j^{\pm 1};p,q)}{\eg(z_j^{\pm2};p,q)}
\frac{dz_1}{z_1}\cdots\frac{dz_n}{z_n}
\\
= \prod_{j=1}^n\Big(\frac{\eg(t^j;p,q)}{\eg(t;p,q)}
\prod_{1\leq m<s\leq 6}\eg(t^{j-1}t_mt_s;p,q )\Big).
\label{SintB}\end{eqnarray}
\end{theorem}
In order to prove formula \eqref{SintB}, consider the following $(2n-1)$-tuple integral
\begin{eqnarray}
&& \makebox[-1em]{}\kappa_{n}\kappa_{n-1}\int_{\T^{2n-1}}
\prod_{1\leq j<k\leq n}\frac{1}{\eg(z_j^{\pm 1} z_k^{\pm 1};p,q)}
 \prod_{j=1}^n\frac{\prod_{r=0}^5\eg(t_rz_j^{\pm 1};p,q)}
{\eg(z_j^{\pm2};p,q)}
\nonumber \\ &&
\times
\prod_{\stackrel{1\leq j\leq n}{1\leq k\leq n-1}}
\eg(t^{1/2}z_j^{\pm 1} w_k^{\pm 1};p,q)
\prod_{1\leq j<k\leq n-1}\frac{1}{\eg(w_j^{\pm 1} w_k^{\pm 1};p,q)}\nonumber \\
&&
\times \prod_{j=1}^{n-1}
\frac{\eg(w_j^{\pm 1} t^{n-3/2}\prod_{s=1}^5t_s;p,q)}
{\eg(w_j^{\pm2},w_j^{\pm 1} t^{2n-3/2}\prod_{s=1}^5t_s;p,q)}
\frac{dw_1}{w_1}\cdots\frac{dw_{n-1}}{w_{n-1}}
\frac{dz_1}{z_1}\cdots\frac{dz_n}{z_n}, \label{compint}
\end{eqnarray}
with the parameters $p, q, t$ and $t_r$, $r=0,\ldots,5,$ lying inside the unit circle
and such that $t^{n-1}\prod_{r=0}^5t_r=pq$. Denote the integral in the left-hand side of
equality \eqref{SintB} by $I_n(t,t_1,\ldots,t_5;p,q)$. Integration over the
variables $w_j$ with the help of formula \eqref{C-typeI} brings expression
\eqref{compint} to the form $\eg^n(t)I_n(t,t_1,\ldots,t_5;p,q)/\eg(t^n)$
(after denoting $t_6=pq/t^{2n-2}\prod_{j=1}^5t_j$).
Because the integrand is bounded on the integration contour, we can
change the order of integrations. As a result, integration over the
variables $z_j$ with the help of formula \eqref{C-typeI} brings
expression \eqref{compint} in the form
$\eg^{n-1}(t) \prod_{0\leq r< s\leq 5} \eg(t_rt_s)
I_{n-1}(t,t^{1/2}t_1,\ldots,$ $t^{1/2}t_5;p,q),$
i.e. we obtain the following recurrence
relation in the dimensionality of the integral of interest $n$:
$$
I_n(t,t_1,\ldots,t_5;p,q)= \frac{\eg(t^n;p,q)}{\eg(t;p,q)}
\prod_{0\leq r<s\leq 5}\makebox[-0.5em]{}\eg(t_rt_s;p,q)\;
I_{n-1}(t,t^{1/2}t_1,\ldots,t^{1/2}t_5;p,q).
$$
Iterating it with known initial condition \eqref{ell-int} for $n=1$,
one obtains formula \eqref{SintB}.

Integral \eqref{SintB} was constructed by van Diejen and the author
in \cite{die-spi:elliptic} and classified as of type II in
\cite{die-spi:selberg} where from the described proof is taken.
This proof models Anderson's derivation of the Selberg integral described
in \cite{aar} (see Theorem 8.1.1 and Sect. 8.4). It also represents a direct
generalization of Gustafson's method \cite{gus:some} of derivation of the
multiple $q$-beta integral obtained from formula \eqref{SintB}
after expressing $t_6$ via other parameters, removing the multipliers $pq$
with the help of the reflection formula for $\eg(z;p,q)$, and taking the limit
$p\to0$. A number of further limits in parameters
leads to the Selberg integral -- one of the most important known integrals because
of many applications in mathematical physics \cite{fw}.
Therefore formula \eqref{SintB} represents an elliptic analogue of
the Selberg integral (an analogous extension of Aomoto's integral described
in Theorem 8.1.2 of \cite{aar} is derived in \cite{rai:trans}).
It can be interpreted also as an elliptic extension
of the $BC_n$ Macdonald-Morris constant term identities.

In analogy with the one dimensional case \cite{spi:theta2}, it is natural
to expect that the multiple elliptic beta integrals define
measures in the biorthogonality relations for some functions of
many variables generalizing relations \eqref{2ib}. In
\cite{rai:trans,rai:abelian}, Rains has constructed a system of
such functions on the basis of integral \eqref{SintB}. These functions
generalize also the Macdonald and Koornwinder orthogonal polynomials,
as well as the interpolating polynomials of Okounkov.  For a related
work see also \cite{cg}.
A systematic investigation of the limiting cases of univariate and multiple
elliptic biorthogonal functions is performed in \cite{br3}.
In this sense, the results obtained in \cite{rai:trans,rai:abelian} represent  to the
present moment the top level achievements of the theory
of elliptic hy\-per\-geo\-met\-ric functions of many variables.
In particular, the following  $BC_n$-generalization of
transformation \eqref{E7-1} was proved in \cite{rai:trans}:
\begin{equation}\label{rains}
I_n(t_1,\ldots,t_8;t;q,p)=I_n(s_1,\ldots,s_8;t;q,p),
\end{equation}
where
\begin{eqnarray*} && \makebox[-2em]{}
I_n(t_1,\ldots,t_8;t;q,p)=\kappa_n\prod_{1\leq j<k\leq 8}\Gamma(t_jt_k;p,q,t)
\\ && \makebox[2em]{} \times
\int_{\T^n}\makebox[-0.5em]{} \prod_{1\le j<k\le n}\!
\frac{\eg(tz_j^{\pm 1}z_k^{\pm 1};p,q)}
{\eg(z_j^{\pm 1}z_k^{\pm 1};p,q)}
\ \prod_{j=1}^n \frac{\prod_{k=1}^8\eg(t_kz_j^{\pm 1};p,q)}
{\eg(z_j^{\pm 2};p,q)}\frac{dz_j}{z_j},
\\ && \makebox[-2em]{}
\left\{
\begin{array}{cl}
s_j =\rho^{-1} t_j,&   j=1,2,3,4  \\
s_j = \rho t_j, &    j=5,6,7,8
\end{array}
\right.;
\quad \rho=\sqrt{\frac{t_1t_2t_3t_4}{pqt^{1-n}}}
=\sqrt{\frac{pqt^{1-n}}{t_5t_6t_7t_8}},
\quad |t|,|t_j|,|s_j|<1,
\end{eqnarray*}
and $\Gamma(z;p,q,t)=\prod_{j,k,l=0}^\infty(1-zt^jp^kq^l)(1-z^{-1}t^{j+1}p^{k+1}
q^{l+1})$ is the elliptic gamma function of the higher level connected
to the Barnes gamma function $\Gamma_4(u;\mathbf{\omega})$.
In \cite{stat}, this symmetry transformation is represented in the
star-star relation form of solvable models of statistical mechanics,
and equality \eqref{SintB} is represented in the star-triangle relation form
which used an elliptic gamma function of even higher order
related to $\Gamma_5(u;\mathbf{\omega})$-function.

There are about 10 proven exact evaluations of elliptic beta integrals  on root systems.
In particular, in \cite{spi:theta2} the author has constructed three
different integrals for the $A_n$ root system (two of them have
different evaluation formulas for even and odd values of $n$).
In \cite{spi-war:inversions}, Warnaar and the
author have found one more $A_n$-integral which appeared to be
new even after degeneration to the $q$- and plain
hy\-per\-geo\-met\-ric levels. Another $BC_n$-integral has been
constructed in \cite{bult:trafo,rai:lit}.
Very many new multiple elliptic beta integrals and symmetry transformations
for their higher order generalizations were conjectured in \cite{SV,SV2}.

Let us describe a generalization of the elliptic beta integral \eqref{SintB}.
Take 10 parameters $p,q,$ $t, s,$ $t_j,$ $s_j$, $j=1,2,3,$
of modulus less than 1 such that
$(ts)^{n-1}\prod_{k=1}^3t_ks_k=pq$ and define the $A_n$-integral
\begin{eqnarray} \lab{AIIa}
&& I_n(t_1,t_2,t_3;s_1,s_2,s_3;t;s;p,q)=\frac{(p;p)_\infty^n(q;q)_\infty^n}
{(n+1)!(2\pi\textup{i})^n}
\\ && \makebox[-1em]{} \times
\int_{\T^n} \prod_{1\leq i<j\leq n+1}
\frac{\eg(tz_iz_j,sz_i^{-1}z_j^{-1};p,q)}
{\eg(z_iz_j^{-1},z_i^{-1}z_j;p,q)}
\prod_{j=1}^{n+1}\prod_{k=1}^3\eg(t_kz_j,s_kz_j^{-1};p,q)
\prod_{j=1}^n \frac{dz_j}{z_j},
\nonumber\end{eqnarray}
where $\prod_{j=1}^{n+1}z_j=1$. Then for odd $n$ one has
\begin{eqnarray} \nonumber
&& \makebox[-2em]{}
I_n(t_1,t_2,t_3;s_1,s_2,s_3;t;s;p,q)=
\eg(t^{\frac{n+1}{2}},s^{\frac{n+1}{2}};p,q)
\\  && \makebox[-1em]{} \times
\prod_{1\leq i<k\leq 3} \eg(t^{\frac{n-1}{2}}t_it_k,
s^{\frac{n-1}{2}}s_is_k;p,q)
\prod_{j=1}^{(n+1)/2}\prod_{i,k=1}^3\eg((ts)^{j-1}t_is_k;p,q)
\lab{AIIa-odd} \\ && \makebox[-1em]{}  \times
\prod_{j=1}^{(n-1)/2}\Big(\eg((ts)^j;p,q)\prod_{1\leq i<k\leq 3}
\eg(t^{j-1}s^jt_it_k,t^js^{j-1}s_is_k;p,q)\Big),
\nonumber\end{eqnarray}
and for even $n$ one has
\begin{eqnarray} \nonumber
&& I_n(t_1,t_2,t_3;s_1,s_2,s_3;t;s;p,q)
= \prod_{i=1}^3\eg(t^{\frac{n}{2}}t_i,s^{\frac{n}{2}}s_i;p,q) \qquad\qquad\qquad
\\  &&\makebox[4em]{} \times
\eg(t^{\frac{n}{2}-1}t_1t_2t_3,s^{\frac{n}{2}-1}s_1s_2s_3;p,q)
\prod_{j=1}^{n/2}\Big(\eg((ts)^j;p,q)
\lab{AIIa-even}\\ && \makebox[4em]{}
\times \prod_{i,k=1}^3\eg((ts)^{j-1}t_is_k;p,q)
\prod_{1\leq i<k\leq 3}\eg(t^{j-1}s^jt_it_k,t^js^{j-1}s_is_k;p,q)\Big).
\nonumber \end{eqnarray}
These $A_n$-elliptic beta integrals were discovered by the author
in \cite{spi:theta2}. As indicated in \cite{SV},
the limit $s\to 1$ reduces the odd $n$ evaluation formula \eqref{AIIa-odd} to
\eqref{SintB}, i.e. we have a generalization of the elliptic
Selberg integral of \cite{die-spi:elliptic,die-spi:selberg}.
The observation that the type II $BC_n$-hypergeometric
identities can be obtained from the type II
relations for $A_{2n-1}$ and $A_{2n}$ root systems was first made
in \cite{spi-war:psi} at the level of multiple $q$-hypergeometric series.
It was also suggested there that the multiple elliptic biorthogonal
rational functions associated with elliptic beta integrals
\eqref{AIIa-odd} and \eqref{AIIa-even},
the existence of which was conjectured by the author long ago
\cite{spi:theta2}, should also generalize the Rains biorthogonal
functions \cite{rai:trans,rai:abelian} to $A_n$ root system.

In  \cite{SV},  the following symmetry transformation was conjectured
for a two parameter extension of the $A_{2n-1}$-integral \eqref{AIIa}:
\begin{eqnarray}\nonumber
&& \makebox[-2em]{}
\int_{\mathbb{T}^{2n-1}} \prod_{1 \leq j < k \leq 2n}
\frac{\Gamma(tz_jz_k,sz_j^{-1}z_k^{-1};p,q)}
{\Gamma(z_j^{-1}z_k,z_jz_k^{-1};p,q)} \prod_{j=1}^{2n}
\prod_{k=1}^{4} \Gamma(t_kz_j,s_kz_j^{-1};p,q)\prod_{j=1}^{2n-1}
\frac{dz_j}{z_j}
 \\\label{SU2N_1}  && \makebox[-2em]{}
=\prod_{1 \leq i < j \leq 4}\Big(
\Gamma(s^{n-1}s_is_j,t^{n-1}t_it_j;p,q)
\prod_{m=0}^{n-2} \Gamma(t (st)^{m}s_is_j,s (st)^{m}t_it_j;p,q)\Big)
\\  && \makebox[-4em]{} \times
\int_{\mathbb{T}^{2n-1}}  \prod_{1 \leq j < k \leq 2n}
\frac{\Gamma(sz_jz_k,tz_j^{-1}z_k^{-1};p,q)}
{\Gamma(z_j^{-1}z_k,z_jz_k^{-1};p,q)}
\prod_{j=1}^{2n}\prod_{k=1}^{4} \Gamma\Big(\sqrt[4]{\frac{S}{T}}
t_kz_j,\sqrt[4]{\frac{T}{S}}s_kz_j^{-1};p,q\Big) \prod_{j=1}^{2n-1}
\frac{dz_j}{z_j},
\nonumber\end{eqnarray}
where $\prod_{j=1}^{2n}z_j=1$, the balancing condition reads
$(st)^{2n-2} ST =(pq)^2$, $S=\prod_{k=1}^4 s_k $ and
$T =\prod_{k=1}^4 t_k,$ and $|s|,|t|,|s_j|,|t_j|,|\sqrt[4]{T/S}s_j|,|\sqrt[4]{S/T}t_j|<1$.
As shown in \cite{SV}, for $s\to 1$ this formula passes
to the Rains transformation \eqref{rains}  and
there are also two more similar symmetry transformations.
Because the integrals in \eqref{SU2N_1} have only $S_4\times S_4 \times S_2$
permutational symmetry in the parameters instead of the $S_8$-group
of \eqref{SintB}, these three Weyl group transformations
lead not to the $E_7$-group, but to a much smaller group.
Consideration of the analogous symmetry transformations for
integrals on the root system $A_{2n}$ has not been completed yet.

\medskip {} {\bf \em An elliptic Fourier transform and a Bailey lemma.}
The Bailey chains, discovered by Andrews, serve
as a powerful tool for building constructive identities for hy\-per\-geo\-met\-ric
series (see Ch. 12 in \cite{aar}). They describe mappings of given sequences
of numbers to other sequences with the help of matrices admitting
explicit inversions.  So, the most general Bailey chain for the univariate
$q$-hy\-per\-geo\-met\-ric series suggested in \cite{and:bailey} is
connected to the matrix built from the $_8\varphi_7$ Jackson sum \cite{bre}.
An elliptic generalization of this chain for the $_{s+1}V_s$-series
was built in \cite{spi:bailey1}, but we do not consider it here, as well as
its complement described in \cite{war:extensions}.
Instead we present a generalization of the  formalism of
Bailey chains to the level of integrals discovered in \cite{spi:bailey2}.

Let us define an integral transformation, which we call an
elliptic Fourier transformation,
\begin{equation}
\beta(w,t)=M(t)_{wz}\alpha(z,t):=\frac{(p;p)_\infty(q;q)_\infty}{4\pi\textup{i}}\int_\mathbb{T}
\frac{\Gamma(tw^{\pm1}z^{\pm1};p,q)}
{\Gamma(t^2,z^{\pm2};p,q)}\alpha(z,t)\frac{dz}{z},
\label{EFT}\end{equation}
where $|tw|,|t/w|<1$ and $\alpha(z,t)$ is an analytical function of variable $z\in\T$.
For convenience we use matrix notation for the $M$-operator and assume in its
action an integration over the repeated indices.
The functions $\alpha(z,t)$ and $\beta(z,t)$ related in the indicated
way are said to form an integral elliptic Bailey pair with respect
to the parameter $t$. An integral analogue of the Bailey lemma, providing an
algorithm to build infinitely many Bailey pairs out of a given one, in this
case has the following form.
\begin{theorem}
Let $\alpha(z,t)$ and $\beta(z,t)$ form an integral elliptic Bailey pair
with respect to the parameter $t$. Then for $|s|,|t|<1, |\sqrt{pq}y^{\pm1}|<|st|$
the functions
\begin{eqnarray} &&
\alpha'(w,st)=D(s;y,w)\alpha(w,t),\quad D(s;y,w)
=\Gamma(\sqrt{pq}s^{-1}y^{\pm1}w^{\pm1};p,q),
\\  &&
\beta'(w,st)=D(t^{-1};y,w) M(s)_{wx}D(st;y,x)\beta(x,t),
\label{BL}\end{eqnarray}
where $w\in \T$, form an integral elliptic Bailey pair with
respect to the parameter $st$.
\end{theorem}

The operators $D$ and $M$ obey nice algebraic properties.
Reflection equation for the elliptic gamma function yields
$D(t^{-1};y, w)D(t;y,w)=1$. As shown in \cite{spi-war:inversions},
under certain restrictions onto the parameters and contours
of integration of the operators $M(t^{-1})_{wz}$ and $M(t)_{wz}$
they become inverses of each other.  Passing to the real integrals
\cite{S5,stat} one can use the generalized functions
and find $M(t^{-1})M(t)=1$ in a symbolic notation where "1" means
a Dirac delta-function. This $t\to t^{-1}$ inversion
resembles the key property of the Fourier
transform and justifies the name ``elliptic Fourier transformation".
The second Bailey lemma given in \cite{spi:bailey2}
is substantially equivalent to this inversion statement.

The conjectural equality $\beta'(w,st)=M(st)_{wz}\alpha'(z,st)$
boils down to the operator identity known as the star-triangle relation
\begin{equation}
M(s)_{wx}D(st;y,x)M(t)_{xz}=D(t;y,w)M(st)_{wz}D(s;y,z),
\label{STR}\end{equation}
which was presented in \cite{spi:umnrev} as a matrix relation (6.5).
After plugging in explicit expressions for $M$ and $D$-operators
one can easily verify \eqref{STR} by using the elliptic beta integral
evaluation formula, which proves the Theorem.

Let us take four parameters $\mathbf{t}=(t_1,t_2,t_3,t_4)$ and consider
elementary transposition operators $s_1,s_2,s_3$ generating
the permutation group $\mathfrak{S}_4$:
$$
s_1(\mathbf{t})=(t_2,t_1,t_3,t_4), \quad
s_2(\mathbf{t})=(t_1,t_3,t_2,t_4), \quad
s_3(\mathbf{t})=(t_1,t_2,t_4,t_3).
$$
Define now three operators $\mathrm{S}_1(\mathbf{t}), \mathrm{S}_2(\mathbf{t})$
and $\mathrm{S}_3(\mathbf{t})$ acting in the space
of functions of two complex variables $f(z_1,z_2)$:
\begin{eqnarray*} && \makebox[-2em]{}
[\mathrm{S}_1(\mathbf{t})f](z_1,z_2):=M(t_1/t_2)_{z_1z}f(z,z_2), \quad
[\mathrm{S}_3(\mathbf{t})f](z_1,z_2):=M(t_3/t_4)_{z_2z}f(z_1,z),
\\ && \makebox[4em]{}
[\mathrm{S}_2(\mathbf{t})f](z_1,z_2):=D(t_2/t_3;z_1,z_2)f(z_1,z_2).
\end{eqnarray*}
As shown in \cite{DS}, these three operators generate the group $\mathfrak{S}_4$,
provided their sequential action is defined via a cocycle condition
$\mathrm{S}_j\mathrm{S}_k:=\mathrm{S}_j(s_k(\mathbf{t}))\mathrm{S}_k(\mathbf{t}).$
Then one can verify that the Coxeter relations
\begin{equation}
\mathrm{S}_j^2=1, \quad \mathrm{S}_i\mathrm{S}_j=\mathrm{S}_j\mathrm{S}_i \
\text{ for } \ |i-j|>1, \quad
\mathrm{S}_j\mathrm{S}_{j+1}\mathrm{S}_j
=\mathrm{S}_{j+1}\mathrm{S}_j\mathrm{S}_{j+1}
\label{coxeter}\end{equation}
are equivalent to the algebraic properties of the Bailey lemma entries,
with the last cubic relation being equivalent to \eqref{STR}.
Thus the Bailey lemma of \cite{spi:bailey1,spi:bailey2}
is equivalent to the Coxeter relations for a permutation group generators \cite{DS}.

The above theorem is used analogously to the Bailey lemma for series \cite{aar}:
one takes initial $\alpha(z,t)$ and $\beta(z,t)$,
found, say, from formula \eqref{ell-int}, and generates
new pairs with the help of the described rules applied to different variables.
Equality \eqref{EFT} for these pairs leads to a tree of
identities for elliptic hy\-per\-geo\-met\-ric integrals of different
multiplicities. As an illustration, we would like to give one nontrivial
relation. With the help of formula \eqref{ell-int}, one can easily
verify the validity of the following recurrence relation
\ba \label{rec-int}
&& I^{(m+1)}(t_1,\ldots,t_{2m+8})=\frac{\prod_{2m+5\leq k<l\leq 2m+8}
\Gamma(t_kt_l;p,q)}{\Gamma(\rho_m^2;p,q)}
\\ && \makebox[1em]{} \times
\kappa\int_\T\frac{\prod_{k=2m+5}^{2m+8}\Gamma(\rho_m^{-1}t_kw^{\pm 1};p,q)}
{\Gamma(w^{\pm2};p,q)}
I^{(m)}(t_1,\ldots,t_{2m+4},\rho_m w,\rho_m w^{-1})\frac{dw}{w},
\nonumber\ea
where $\rho_m^2=\prod_{k=2m+5}^{2m+8}t_k/pq$ and the integral
$I^{(m)}$ was defined in \eqref{gen-int}. By an appropriate change of notation,
one obtains a concrete realization of the Bailey pairs:
$\alpha\propto I^{(m)}$ and $\beta\propto I^{(m+1)}$.
For $m=0$, substitution of the explicit expression \eqref{ell-int} for
$I^{(0)}$ in the right-hand side of \eqref{rec-int}
yields identity  \eqref{E7-1}. Other interesting consequences of
the recursion \eqref{rec-int} (an elliptic analogue of formula
(2.2.2) in \cite{aar}) are considered in \cite{spi:thesis,S5}.
Various generalizations of the elliptic Fourier transformation \eqref{EFT}
to root systems and their inversions are described in \cite{spi-war:inversions}.

\medskip {} {\bf \em Connection to the representation theory.}
Plain hypergeometric functions are connected to matrix
elements of the representations of standard Lie groups
(see, e.g., Sect. 9.14 in \cite{aar} where the Jacobi polynomials case
is considered).
Some of the $q$-special functions have been interpreted in a
similar way in connection to quantum groups. Therefore it is
natural to try to construct elliptic hypergeometric functions
from the representations of ``elliptic quantum groups".
The current top result along these lines
was obtained in \cite{ros:algebra}, where the terminating elliptic
hypergeometric series of type I on the $A_n$ root system was
constructed as matrix elements for intertwiners between
corepresentations of an elliptic quantum group. However, the whole
construction is quite complicated and the elliptic hypergeometric
integrals have not been treated in this way yet.

A qualitatively new group-theoretical interpretation of
the elliptic hypergeometric functions has emerged, again,
from mathematical physics (see \cite{DO,gprr,SV,SV2}
and references therein). It directly connects the elliptic hypergeometric
integrals to the representations of standard Lie groups.
Take a Lie group $G\times F$ and a set of its irreducible
representations including the distinguished
representation $\adj_G$, adjoint for group $G$ and trivial for $F$
(the ``vector" representation).
Consider the following function of this group characters:
\begin{equation}
I(y;p,q) \ = \ \int_{G} d \mu(z)\,
\exp \Big( \sum_{n=1}^{\infty}
\frac 1n \ind\big(p^n ,q^n, z^n , y^ n\big ) \Big),
\label{int_G}\end{equation}
where $d \mu(z)$ is the $G$-group invariant (Haar) measure and
\begin{eqnarray}\nonumber  &&
\ind(p,q,z,y) =  \frac{2pq - p - q}{(1-p)(1-q)} \chi_{\adj_G}(z)
 \\ && \makebox[2em]{}
+ \sum_j \frac{(pq)^{r_j}\chi_{R_F,j}(y)\chi_{R_G,j}(z) - (pq)^{1-r_j}
\chi_{{\bar R}_F,j}(y)\chi_{{\bar R}_G,j}(z)}{(1-p)(1-q)}
\label{1ind} \end{eqnarray}
with some fractional numbers $r_j$. Here $\chi_{\adj_G}(z)$ and $\chi_{R_G,j}(z)$,
$\chi_{R_F,j}(y)$ are the characters of the ``vector" and all other (``chiral")
representations, respectively. They depend
on the maximal torus variables $z_a$, $a=1,\ldots,\text{rank}\, G$,
and $y_k$, $k=1,\dots,\text{rank}\, F$.

For  $G=SU(N)$ one has $z=(z_1,\ldots,z_N),$ $\prod_{j=1}^Nz_j=1$, and
\begin{eqnarray*} &&
\int_{SU(N)} d\mu(z) \ = \   \frac{1}{N!} \int_{\mathbb{T}^{N-1}}
\Delta(z) \Delta(z^{-1}) \prod_{a=1}^{N-1} \frac{dz_a}{2 \pi \textup{i} z_a},
\end{eqnarray*}
where $\Delta(z) \ = \ \prod_{1 \leq a < b \leq N} (z_a-z_b)$, and
$\chi_{SU(N),\adj}(z)= (\sum_{i=1}^N z_i)(\sum_{j=1}^N z_j^{-1})-1$.

For special sets of representations entering the sum $\sum_j$ in \eqref{1ind}
and some fractional numbers $r_j$ formula \eqref{int_G} yields {\em all} known elliptic
hypergeometric integrals with interesting properties. It has even deeper
group-theoretical meaning in the context of the representation theory
of superconformal group $SU(2,2|1)$, where $2r_j$ coincide with the
eigenvalues of $U(1)_R$-subgroup generator (``$R$-charges") and $p,q$
are interpreted as group parameters for generators commuting with
a distinguished pair of supercharges (see the next section).

Take the elliptic beta integral \eqref{ell-int} and rewrite it as $I_{\text{lhs}}=I_{\text{rhs}}$,
where $t_k=(pq)^{1/6}y_k$, $k=1,\ldots,6$. Then $I_{\text{lhs}}$ is obtained from \eqref{int_G}
for $G=SU(2),$ $F=SU(6)$ with two representations: the ``vector" one
$(\adj, 1)$ with $\chi_{SU(2),\adj}(z)=z^2+z^{-2}+1$ and the fundamental one
$(f, f)$ with $\chi_{SU(2),f}(z)=z+z^{-1},$ $r_f=1/6,$ and
$$
\chi_{SU(6),f}(y)=\sum_{k=1}^6y_k,
\quad \chi_{SU(6),\bar f}(y)=\sum_{k=1}^6y_k^{-1},
\quad \prod_{k=1}^6y_k=1.
$$
The latter constraint on $y_k$ is nothing else than the balancing condition
for the integral in appropriate normalization of parameters, i.e. this notorious
condition is equivalent to the demand that
the determinant of special unitary matrices is equal to 1.
For $I_{\text{rhs}}$ one has $G=1,$ $F=SU(6)$ with single representation
$T_A:\, \Phi_{ij}=-\Phi_{ji},\ i,j=1,\ldots,6,$ with
$$
\chi_{SU(6),T_A}(y)=\sum_{1\leq i<j\leq 6}y_iy_j,
\qquad r_{T_A}=1/3.
$$
The elliptic beta integral evaluation formula thus proves the equality
of two character functions on
different groups with different sets of representations. All known analogous
relations between integrals can be interpreted in this way.
Since the elliptic hypergeometric integrals are expected to define
automorphic functions in the cohomology class of the group $SL(3,\Z)$,
this could mean the equivalence of two differently defined
automorphic functions, which is a new type of
group-theoretical duality. A physical interpretation
of this construction is described in the next section.

\medskip {} {\bf \em Applications in mathematical physics.}
The most important known physical application of elliptic hypergeometric
integrals has been found in four dimensional supersymmetric quantum field
theories, where they emerge as superconformal indices.

For $\mathcal{N}=1$ supersymmetric theories the full symmetry group is
$G_{\text{full}}=SU(2,2|1)\times G\times F$, where the space-time symmetry group
is generated by $J_i, \overline{J}_i$, $i=1,2,3$ ($SU(2)$ subgroup generators,
or $SO(3,1)$-group Lorentz rotations), $P_\mu, Q_{\alpha},
\overline{Q}_{\dot\alpha}$, $\mu=0,\ldots,3$, $\alpha, \dot\alpha=1,2$
(supertranslations), $K_\mu, S_{\alpha},\overline{S}_{\dot\alpha}$
(special superconformal transformations),
$H$ (dilations), and $R$ ($U(1)_R$-rotations);
$G$ is a local gauge invariance group and $F$
is a global flavor symmetry group.
The whole set of commutation relations between these operators can
be found, e.g., in \cite{SV}. Choosing a particular pair of
supercharges, say, $Q=\overline{Q}_{1 }$ and $Q^{\dag}=-{\overline S}_{1}$,
one obtains
\begin{equation}
Q Q^{\dag}+Q^{\dag}Q = 2{\mathcal H},\quad Q^2= (Q^{\dag})^2=0,\qquad
\mathcal{H}=H-2\overline{J}_3-3R/2.
\label{susy}\end{equation}
Then the superconformal index (SCI) is defined by the following trace:
\begin{eqnarray}
I(y;p,q) = \text{Tr} \Big( (-1)^{\mathcal F}
p^{\mathcal{R}/2+J_3}q^{\mathcal{R}/2-J_3}
\prod_k y_k^{F_k} e^{-\beta {\mathcal H}}\Big),
\quad \mathcal{R}= H-R/2,
\label{Ind}\end{eqnarray}
where $\mathcal F$ is the fermion number operator ($(-1)^{\mathcal F}$ is
simply a $\Z_2$-grading operator in $SU(2,2|1)$), $F_k$ are the maximal
torus generators of the group $F$, and $p,q,y_k,\beta$ are group parameters.
The trace in \eqref{Ind} is taken over the Hilbert (Fock) space of quantum fields
forming irreducible representations of the group $G_{\text{full}}$.
Because operators $\mathcal{R}, J_3, F_k, {\mathcal H}$ used in the definition of
SCI commute with each other and with $Q, Q^\dag$, non-zero contributions to the trace
may come only from the space of zero modes of the operator $\mathcal{H}$ (or the cohomology
space of $Q$ and $Q^\dag$ operators). Therefore there is no $\beta$-dependence.
Computation of this trace leads to integral \eqref{int_G},
where the integration over $G$ reflects the gauge invariance of SCI.
Function \eqref{1ind} is called the one-particle states index.

Some of the supersymmetric field theories are related to one another by
the Seiberg electric-magnetic dualities \cite{seiberg}, which are not proven
yet despite of many convincing arguments.
Equality of SCIs for such theories was
conjectured by R\"omelsberger and proved in some cases by Dolan and Osborn \cite{DO}
by identifying SCIs with the elliptic hypergeometric integrals.
A related application to topological quantum field theories
(which is using an elliptic hypergeometric integral identity of \cite{bultF4})
is discovered in \cite{gprr}.
In \cite{SV,SV2} many new $\mathcal{N}=1$ supersymmetric dualities have been found
and very many new integral identities
have been conjectured, among which there are relations of a qualitatively new type
(e.g., they involve higher order generalizations of integral \eqref{SintB}
with $t=(pq)^{1/K}$, $K=2,3,\ldots$).

We leave it as an exercise to determine what kind of transformation of
elliptic hypergeometric integrals is hidden behind the equality of SCIs
for the original Seiberg duality \cite{seiberg}. In this
case one has two theories
with $F=SU(M)_l\times SU(M)_r\times U(1)$ (here $U(1)$ is the baryon number
preserving symmetry) and different gauge groups
and representations. The ``electric" theory has the group $G=SU(N)$ and the set
of representations described in the table below:
\begin{center}\begin{tabular}{|c|c|c|c|c|}
\hline
$SU(N)$ & $SU(M)_l$ & $SU(M)_r$ & $U(1)$ & $U(1)_R$ \\
\hline
 $f$ & $f$ & 1 & 1 & $\tilde{N}/M$ \\
 $\overline{f}$ & 1 & $\overline{f}$ & $-1$ & $\tilde{N}/M$ \\
 ${\rm \adj}$  & $1$   &  $1$ &  $0$   &  $1$ \\
\hline
\end{tabular}
\end{center}
where $\tilde{N}=M-N$. The ``magnetic" theory has the group $G=SU(\tilde{N})$
with the representations described in the following table:
\begin{center}
\begin{tabular}{|c|c|c|c|c|}
\hline
 $SU(\tilde{N})$ & $SU(M)_l$ & $SU(M)_r$ & $U(1)$ & $U(1)_R$ \\
\hline
 $f$ & $\overline{f}$ & 1 & $N/\tilde{N}$ & $N/M$ \\
 $\overline{f}$ & 1 & $f$ & $-N/\tilde{N}$ & $N/M$
\\
 1 & $f$ & $\overline{f}$ & 0 & $2\tilde{N}/M$
\\
${\rm \adj}$  & $1$   &  $1$ &  $0$   &  $1$ \\
\hline
\end{tabular}
\end{center}
The last columns of these tables contain the numbers $2r_j$ -- eigenvalues
of the generator of $U(1)_R$-group $R$. The last
rows correspond to the vector superfield representation,
other rows describe chiral superfields.
For $N=2, M=3$ equality
of SCIs is equivalent to the elliptic beta integral, as described in the
previous section. For arbitrary $N$ and $M$ SCIs were computed in \cite{DO}
(see also \cite{SV}).
Physically, the exact computability of SCIs describes
a principally important physical phenomenon -- the confinement of colored particles
in supersymmetric theories of strong interactions. The equality of SCIs
provides presently the most rigorous mathematical justification of the
Seiberg dualities.

In \cite{spi-zhe:spectral}, a discrete integrable system
generalizing the discrete-time Toda chain has been constructed.
A particular elliptic solution of this nonlinear chain equations
has lead to the terminating ${}_{12}V_{11}$-series as a solution
of the Lax pair equations. Derivation of this function from a similarity reduction of
an integrable system equations reflects the es\-sence of a powerful heuristic
approach to all special functions of one variable (it was described in detail in
\cite{spi:thesis} on the basis of a number of other
new special functions constructed in this way).
In \cite{kmnoy}, it was shown that the same ${}_{12}V_{11}$-series appears as a particular
solution of the elliptic Painlev\'e equation discovered by Sakai \cite{sak}.
An analogous role is played by the general solution of the elliptic
hy\-per\-geo\-met\-ric equation \cite{spi:thesis,spi:cs} and some
multiple elliptic hy\-per\-geo\-met\-ric integrals \cite{rai:rec,rai:painleve}.
In \cite{BS}, a different discrete integrable system was deduced from
the semiclassical analysis of the elliptic beta integral.

The first physical interpretation of elliptic hypergeometric integrals was
found in \cite{spi:thesis,spi:cs}, where it was shown that some of
the $BC_n$-integrals describe either special wave functions or normalizations
of wave functions in the Calogero-Sutherland type many body quantum
mechanical models. One can consider in an analogous way the root system $A_n$.
It is natural to expect that all superconformal indices are associated
with such integrable systems \cite{SV}.

Another rich field of applications of elliptic hypergeometric functions
is connected with the exactly solvable models in statistical
mechanics. As mentioned in the introduction, elliptic hypergeometric
series showed up for the first time as solutions of the Yang-Baxter equation
of IRF (interaction round the face) type. The vertex form of the Yang-Baxter
equation naturally leads to the Sklyanin algebra \cite{skl}.
Connection of the elliptic hypergeometric functions
with this algebra is considered in
\cite{DKK,konno,rai:abelian,ros:elementary,ros:sklyanin,S5}. In particular,
in \cite{ros:sklyanin} Rosengren proved an old conjecture of Sklyanin
on the reproducing kernel, and in \cite{S5} an elliptic generalization
of the Faddeev modular double \cite{fad:mod} was constructed.

Let us briefly describe how the $V$-function emerges in this context.
The general linear combination of four Sklyanin algebra generators
can be represented in the form \cite{rai:abelian}
\begin{eqnarray}\nonumber
&& \Delta(\underline{a})
=\frac{\prod_{j=1}^4\theta_1(a_j+u)}
{\theta_1(2u)}e^{\eta\partial_u}+\frac{\prod_{j=1}^4\theta_1(a_j-u)}
{\theta_1(-2u)}e^{-\eta\partial_u},
\nonumber\end{eqnarray}
where $a_j$ and $\eta$ are arbitrary parameters and $e^{\pm\eta\partial_u}f(u)
=f(u\pm\eta)$. The Casimir operators
take arbitrary continuous values, i.e. one deals in general with the continuous
spin representations. The generalized eigenvalue problem
$$
\Delta(a,b,c,d)f(u;\lambda,a,b;s)=\lambda \Delta(a,b,c',d')f(u;\lambda,a,b;s),
$$
where $s=a+b+c+d$ and $c+d=c'+d'$, is exactly solvable and $f$ is given
by a product of 8 elliptic gamma functions. Take the scalar product
$$
\langle f,g\rangle =\kappa \int_{\T}\frac{f(u)g(u)}
{\Gamma(z^{\pm 2};p,q)}
\frac{dz}{z},\qquad z=e^{2\pi i u},\quad p=e^{2\pi i\tau}, \quad q=e^{4\pi i\eta},
$$
and consider the conjugated generalized eigenvalue problem induced
by it
$$
\Delta^*(a,b,c,d)g(u;\mu,a,b;s)=\mu \Delta^*(a,b,c',d')g(u;\mu,a,b;s),
$$
where $\Delta^*$ is defined from the equality
$\langle \Delta f, g\rangle=\langle f, \Delta^* g\rangle$.
Then, the overlap of two dual bases with the same  $s$-parameter
$\langle f, g\rangle$ is equal to $V(\underline{t})$ for appropriately
chosen parameters $t_j$ \cite{S5}.

The general (rank one) solution of the Yang-Baxter equation was derived in \cite{DS}
in the form of an integral operator acting in the space of functions of two complex
variables. The general construction of \cite{DKK}, algebraic properties of the
operators $\mathrm{S}_k$ \eqref{coxeter} (including the fact that the operators
$\mathrm{S}_{1,3}$ are intertwining operators for the Sklyanin algebra generators)
and the elliptic modular double of \cite{S5} played a crucial role in this result.
The Bailey lemma operator identity
\eqref{STR} (or the cubic Coxeter relation) is equivalent to the
star-triangle relation for specific Boltzmann weights considered in \cite{BS}.
In \cite{stat}, the most general solution of the star-triangle
relation connected to the hyperbolic beta integrals was described. Also, it was shown that
the symmetry transformations for elliptic hypergeometric
integrals can be rewritten as the ``star-star" relation (an IRF type
Yang-Baxter equation) leading to new checkerboard type solvable models
of statistical mechanics. The multicomponent generalizations of these models
was proposed there as well. A general relation between different models is
established via the vertex-face correspondence for $R$-matrices.
In this way one gets new Ising-type solvable models
for the continuous spin systems (i.e., two-dimensional quantum field
theories) unifying many previously known examples. The free energy per edge
for the elliptic beta integral model was computed in \cite{BS}. According to the $4d/2d$
correspondence described in \cite{stat} the Seiberg-type dualities
for superconformal indices of
$4d$ supersymmetric gauge field theories are equivalent to the Kramers-Wannier
type duality transformations for elementary cell partition functions,
with the full $2d$ lattice partition functions being equal to superconformal
indices of certain quiver gauge theories.

As a final example, we mention an interesting application of the discrete
elliptic biorthogonal rational functions of \cite{spi-zhe:spectral} to random
point processes related to the statistics of lozenge tilings
of a hexagon (or plane partitions) described in \cite{bgr}.

\medskip

{} {\bf \em Conclusion.}
The main part of the theory of plain hy\-per\-geo\-met\-ric
functions has found a natural elliptic generalization, although
the similarities start to show up for a rather large number
of free parameters and structural restrictions.
We would like to finish by listing some
other achievements of the theory of elliptic hy\-per\-geo\-met\-ric functions.
Multiple elliptic hy\-per\-geo\-met\-ric series were considered for the first time by
Warnaar \cite{war:summation}. We described mostly properties of
the elliptic hy\-per\-geo\-met\-ric integrals, since many results for the
series represent their particular limiting cases being derivable via residue calculus.
A combinatorial proof
of the Frenkel-Turaev summation formula is given in \cite{e-comb}.
Various generalizations of this sum to the root systems were found in
\cite{die-spi:elliptic,ros:elliptic,spi:theta2,war:summation}.
Multivariable analogues of the elliptic Bailey transformation for series
were described in \cite{cg,kaj-nou,rai:trans,ros:elliptic,war:summation}.
Expansions in partial fractions of the ratios of theta functions
and the identities connected to them were considered by Rosengren in
\cite{ros:elliptic} (see also \cite{die-spi:selberg,rai:trans,RS}).
Such expansions played an important role in the method of proving
elliptic beta integral evaluations considered in \cite{spi:short}. Some general
properties of elliptic hypergeometric terms were described in \cite{S4}.
In \cite{spi:6d} a modified elliptic gamma function of the second order was built
and another form of the identity \eqref{rains} having an interesting interpretation
in six-dimensional supersymmetric field theory was found.

The terminating continued fraction generated by the three term recurrence
relation \eqref{ttr} and the Racah type termination condition
was computed in \cite{spi-zhe:rims}. The raising and lowering operators
connected with rational functions were discussed in
\cite{magnus1,magnus2,rai:trans,spi-zhe:grids}.
In particular, in \cite{spi-zhe:grids} it was shown that the
general lowering operator of the first order can exist only
for the elliptic grids.
A systematic investigation of the elliptic determinant formulas
connected to the root systems is performed in the work of Rosengren
and Schlosser \cite{ros-sch3}. Determinants of elliptic
hypergeometric integrals were considered in \cite{RS,spi:theta2}.
Elliptic Littlewood identities were discussed in \cite{rai:lit}
where many quadratic transformations for multiple elliptic hypergeometric
functions were conjectured. Some of these conjectures were proved by van de Bult
\cite{bult:quadr} (quadratic transformations for the univariate
series were derived in \cite{spi:bailey1,war:extensions}). In \cite{zhe:dense},
the elliptic hypergeometric series $_3E_2$ with arbitrary power counting
argument was shown to describe some polynomials with a dense point
spectrum. Connections to the Pad\'e interpolation were analyzed in
\cite{magnus1,spi-zhe:grids,zhe:pad}.

Solutions of various finite
difference equations on the elliptic grid were considered by Magnus
in \cite{magnus1,magnus2}. As shown in \cite{SV2}, reduction of
elliptic hypergeometric integrals to the hyperbolic level leads to
the state integrals for knots in three-dimensional space.
The page limits of the present complement
do not allow the author to cite a number of other interesting results,
an essentially more complete review of the literature
is given in papers \cite{spi:thesis,spi:umnrev} and \cite{SV,SV2}.

To conclude, elliptic hypergeometric functions are universal functions
with important applications in various fields of mathematics
and theoretical physics.
They unify special functions of elliptic and hypergeometric types
under one roof and make them firm, unique, undeformable objects
living in the Platonic world of ideal bodies.

Despite of the very big progress in the development of the theory of
elliptic hypergeometric functions, many open problems still remain.
They include the proof of tens of existing
conjectures on evaluations or symmetry transformations for
integrals, a rigorous definition of
infinite elliptic hy\-per\-geo\-met\-ric series, detailed investigation
of the specific properties of functions when bases are related to
roots of unity,
computation of the nonterminating elliptic hy\-per\-geo\-met\-ric continued
fraction, detailed analysis of the non-self-dual biorthogonal functions
of \cite{spi-zhe:spectral} (still representing the most complicated known
univariate special function of such a kind) and construction
of their multivariable analogues, search of the number theoretic applications
of these functions analogous to those considered in \cite{zudilin},
investigation of their automorphic properties,
higher genus Riemann surface generalizations, and so on.

\smallskip
I am indebted to Yu.~A.~Neretin for the suggestion to write this
complement and to G.~E.~Andrews, R. Askey, and R. Roy for an
enthusiastic support of this idea. I am grateful also to H. Rosengren,
G. S. Vartanov and S. O. Warnaar for useful remarks.
This work is supported in part by the RFBR grants 05-01-01086
and 11-01-00980 (joint with NRU-HSE grant 11-09-0038)
and the Max Planck Institute for mathematics (Bonn).

\end{document}